\documentclass[graybox]{svmult}

\usepackage{type1cm}        
\usepackage{makeidx}         
\usepackage{graphicx}        
\usepackage{multicol}        
\usepackage[bottom]{footmisc}
\usepackage{verbatim}

\usepackage{newtxtext}       %
\usepackage[varvw]{newtxmath}       
\usepackage{overpic}
\usepackage{tikz}
\usepackage{amscd,latexsym,subfigure,verbatim,epsfig,supertabular}
\usepackage[graph,frame,poly,arc]{xy}
\usepackage{calc}
\usepackage{amsmath}
\usepackage{graphicx}

\makeindex             

\newtheorem{defn0}{Definition}[subsection]

\theoremstyle{definition}
\newtheorem{example0}[definition]{Example}

\newcommand{\NN}{{\mathbb N}}

\newcommand{\ZZ}{{\mathbb Z}}
\newcommand{\RR}{{\mathbb R}}

\newcommand{\pp}{{\mathbb P}}

 \newcommand{\cc}{\mathcal C}
\newcommand{\Sr}{S^r(\Delta)} 
\newcommand{\Srd}{S^r_d(\Delta)} 
\newcommand{\Srhat}{S^r(\hat \Delta)} 
 
\newcommand{\stv}{\mathop{\rm st(v)}\nolimits} 
 
\newcommand{\im}{\mathop{\rm im}\nolimits}

\newcommand{\coker}{\mathop{\rm coker}\nolimits}
\newcommand{\AS}{\mathop{\rm AS}\nolimits}

\begin{document}
\title*{The algebra of splines: duality, group actions and homology}

\author{Martina Lanini, Hal Schenck, Julianna Tymoczko}
\institute{Martini Lanini \at Dipartimento di Matematica, Universita di Roma Tor Vergata,  \email{lanini@mat.uniroma2.it} \and Hal Schenck \at Department of Mathematics, Auburn University, \email{hks0015@auburn.edu} \and Julianna Tymoczko \at Department of Mathematics, Smith College, \email{jtymoczko@smith.edu}}
\maketitle


\begin{abstract}
\noindent This survey gives an overview of three central algebraic themes related to the study of splines: duality, group actions, and homology. Splines are piecewise polynomial functions of a prescribed order of smoothness on some subdivided domain $D \subseteq \RR^k$, and appear in applications ranging from approximation theory to geometric modeling to numerical analysis. Alternatively, splines can be interpreted as a collection of polynomials labeling the vertices of a (combinatorial) graph, with adjacent vertex-labels differing by a power of an affine linear form attached to the edge.  In most cases of interest, the subdivided domain is essentially dual to the combinatorial graph, and these two characterizations of splines coincide.   

Properties of splines depend on combinatorics, topology, geometry, and symmetry of a simplicial or polyhedral subdivision $\Delta$ of a region $D \subseteq \RR^k$, and are often quite subtle. We describe how duality, group actions, and homology --- techniques which play a central role in many areas of both pure and applied mathematics -- can be used to illuminate different questions about splines.

Our target audience is nonspecialists:  
we provide a concrete introduction to these methods, and illustrate them with 
many examples in the context of splines.  We also provide a tutorial on computational aspects: all of the objects appearing in this note may be studied using open source computer algebra software.
\end{abstract}

\section{Introduction}\label{sec:intro}
\vskip .1in
Multivariate splines are important tools in numerical analysis and approximation theory. Despite an extensive literature on the subject, numerous open questions remain  about both theoretical and practical aspects, ranging from finding their dimension to constructing local bases and determining their approximation power.  

The core problem is to construct
finite-dimensional spaces of functions that can approximate
complicated or unknown functions well.  Such spaces are especially important
for scientific computing, where they are used in computer-aided geometric design,
data fitting, solving partial differential equations by the
finite-element method, and the isogeometric analysis paradigm. Historically, polynomials provided a first solution to this problem; spaces of splines are a refinement of this solution that are much more efficient and effective. 

An important property of a spline space is its dimension as a vector space, which can be subtle to determine. Some spline spaces have {\em stable} dimension, in the sense that the dimension only depends on the degree $d$, the order of smoothness $r$, and the combinatorial --- or other easy-to-check --- properties of the subdivision $\Delta$, but not on e.g. precisely how $\Delta$ is embedded in Euclidean space. Many useful spline spaces are not stable. The next priority for a useful spline space is identifying optimal approximation power.  In addition, we'd like to be able to construct any function on the spline space locally on each of the elements of the partition. Optimal approximation power turns out to be related \cite{ls} to whether it is possible to construct stable bases with local support for the subdivision $\Delta$.

The interplay between local and global constructions, the two parameters of degree and smoothness, and the effects of varying the subdivision $\Delta$ are at the heart of questions about multivariate splines.  They also connect directly to algebraic geometry, representation theory, and other fields. 

For instance, splines of low degree and high order of smoothness tend to be the most useful in 
applications --- e.g., solving interpolation and approximation problems in $\RR^k$ is much easier with spline spaces of low degree.  
Results from both approximation theory and algebraic geometry show that subdivisions $\Delta$ with many symmetries often have more low-degree splines.  Symmetry in turn implies that a group action is lurking in the background, suggesting a natural connection to representation theory, as discussed in \S\ref{sec:groupActions}. 

 This brings us to another important open question with immediate applications in PDEs, as well as other fields: finding good subdivisions $\Delta$ of a region.  Subdividing an initial subdivision $\Delta$ creates more symmetries (and low-degree splines), but comes with the cost of solving more complicated equations. 
We discuss a number of variants on this theme, including splines on polyhedral or semialgebraic subdivisions, supersmoothness, high dimension, and special partitions such as crosscut or root systems.

\subsection{Historical background} 
Much of the foundational work on splines was developed by numerical analysts using classical methods, in particular Bernstein and B\'ezier techniques. Manni and Sorokina \cite{carlaTanya} provide a survey of these techniques elsewhere in this volume. 

A conjecture of Strang \cite{strang} led Billera \cite{bTAMS} to introduce methods of homological algebra to the study of splines; Billera's work on the dimension question for planar splines of smoothness order one was recognized with the 1994 Fulkerson prize. An introduction to homological methods for non-specialists is the content of \S2 of this paper. 

Splines are also of keen interest to researchers in geometry, topology, and representation theory, where the work of Goresky, Kottwitz and MacPherson on spaces with a nice group action shows that they arise as equivariant cohomology rings.  We sketch this history here and describe it in more detail in \S 3.

In a seminal 1998 paper \cite{GKM98}, Goresky, Kottwitz and MacPherson proved that, under certain technical conditions, the equivariant cohomology of a variety $X$ acted upon by a torus $T$ depends only on a graph $\mathcal{G}_{(X,T)}$ whose vertices are zero-dimensional orbits of $T$ and whose edges are one-dimensional orbits of $T$. 

 Goresky, Kottwitz and MacPherson show that the equivariant cohomology may be identified with a certain set of splines $C^{0}_\infty(\mathcal{G}_{(X,T)})$ of polynomials of unbounded degree and 0-smoothness gluing conditions. The {\em Alfeld split} $\AS^n$ of an $n-$simplex is obtained by taking a barycentric subdivision of a simplex, or what an algebraist would describe as placing a single interior vertex inside a simplex and then coning with the boundary.  We recover the splines $C^{0}_\infty(\AS^n)$ for the $n$-th Alfeld split if $X=\mathbb{P}_{\mathbb{C}}^n$ (see \cite{Ty16}). This observation led to the definition of generalized splines on a graph \cite{GTV16}, which we believe deserve further investigation. 
 
 The GKM description of equivariant cohomology has been especially fruitful in the field of Schubert calculus, which studies the ring structure of the cohomology for varieties of representation theoretic origin (especially flag varieties, using a basis induced by the Schubert varieties contained within them).  Varieties studied in Schubert calculus have a ``nice" cohomology basis, in the sense that the basis generates $C^{0}_\infty(\mathcal{G}_{(X,T)})$ as a free module over an appropriate polynomial ring $R$; this provides a good control on the cohomology ring. The quest for nice bases for more general varieties has interested many researchers \cite{Go14, HTy17, LP20}. Moreover, GKM techniques led to the realisation of symmetric group representations on cohomology modules \cite{LP21, Ty08}. 
 
 An exciting research direction is to investigate which of these techniques and results extend to the setting of degree-bounded polynomials and arbitrary smoothness gluing conditions, that is to $C^{r}_k(\mathcal{G}_{(X,T)})$. For example, the symmetric group actions on $C^{0}_\infty(\mathcal{G}_{(X,T)})$ from \cite{LP21} and \cite{Ty08} preserve the subspace $C^{r}_k(\mathcal{G}_{(X,T)})$ --- so what is this subrepresentation? A related question asks: does the Schubert basis of $C^{0}_\infty(\mathcal{G}_{(X,T)})$ as a free module over $R$ provide information about the dimension of  $C^{r}_k(\mathcal{G}_{(X,T)})$? The answer is yes when $r=0$, but what about general $r$? We touch on these questions more in \S 3.

\subsection{Definitions and Running Examples}
With history behind us, we now introduce our cast of characters and running examples. Our goal is to strike a balance between framing the problem in general while at the same time giving concrete examples that bring the players to life. In short, definitions are general, while examples are specific.
\begin{definition}\label{TheSubdivision}
Let $\Delta$ be a simplicial, polyhedral, or semi-algebraic complex embedded in $\RR^k$, and $D$ the subspace of $\RR^k$ consisting of the union of the cells of $\Delta$.
\end{definition}

On the most abstract level, $\Delta$ could be a CW-complex (see for example \S 5.3 of \cite{sBook2}).  In practice, this level of complexity is overkill. Our primary focus in this survey is on simplicial complexes but we touch on polyhedral and semialgebraic complexes in \S 4, which focuses on open problems. 

We give the formal definition of a simplicial complex in \S 2.2.2 though we often don't even use that level of generality: in most common applications, including Examples~\ref{ex:first} and~\ref{ex:second} below, $\Delta$ is a triangulation of a region in the plane or a tetrahedral decomposition of a region in $\RR^3$.  

\begin{definition}\label{SplineDef}
Let $\Delta$ be a simplicial, polyhedral, or semi-algebraic complex embedded in $\RR^k$. A {\em spline of smoothness} $r \in \NN$ on $\Delta$ is a function $F$ on the subspace $D$ of Definition~\ref{TheSubdivision} which is $r$-differentiable on the interior of $D$, and such that on each maximal cell $\sigma$ the restriction $F|_\sigma$ is given by a polynomial.
\end{definition}

Splines are also referred to as piecewise polynomial functions on $\Delta$.  The set of splines of smoothness $r \in \NN$ such that each $F|_\sigma$ has degree at most $d$ is a vector space denoted by $\Srd$. We assemble these vector spaces into one large set as follows:
\begin{equation}\label{SplineModule}
\Sr = \bigcup\limits_{d=0}^{\infty} S^r_d(\Delta). 
\end{equation}

In fact, $\Sr$ is more than just the union of its parts.  Let $R=\RR[x_1,\ldots,x_k]$.  Note that for every subdivision $\Delta \subseteq \RR^k$, the splines $\Sr$ contain a copy of $R$ that is constructed by assigning the same polynomial $P \in R$ to every face $\sigma_i \in \Delta$. We make this remark because with this understanding of $P \in R$ we have:
\begin{equation}\label{isModule}
\begin{array}{ccc}
F_1\mbox{ and }F_2 \in \Sr & \longrightarrow & F_1+F_2 \in \Sr.\\
P \in R \mbox{ and } F \in \Sr & \longrightarrow & P \cdot F \in \Sr.
\end{array}
\end{equation}
These two equations mean that $\Sr$ is a {\em module} over $R$. This in turn means we can use algebraic tools to shed light on the structure of $\Sr$, as we'll see in \S 2.

Now that we've had a solid dose of terminology, it's time to bring things to life with some examples.
\pagebreak

\begin{example0}\label{ex:first}
Consider the three triangulations below. 
From a combinatorial perspective, $\Delta_1, \Delta_2$, and $\Delta_3$ are all the same: each simplicial complex consists of four triangles and four interior edges. But we consider these simplicial complexes {\em along with their embedding in} $\RR^2$, and the geometry at the central vertex clearly differs. 

\begin{figure}[ht]\label{3planarTriangulations}
\vskip -.1in
  \hskip -1.1in \includegraphics[width=7in]{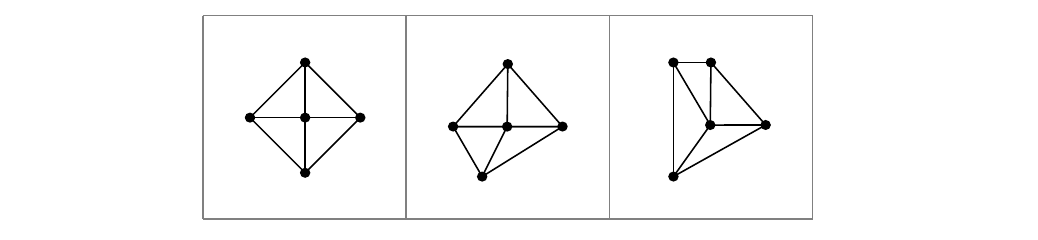}
\begin{center}\caption{Three different geometries $\Delta_1, \Delta_2, \Delta_3$ on the star of 4 vertices} \label{figure: star triangulations}
\end{center}
\end{figure}
\vskip -.3in
Without loss of generality we may assume the central vertex is located at the origin; we call the geometry at a vertex \textit{local geometry}.  The difference in the local geometry of the three triangulations is reflected in the dimensions of the spline spaces, which we record in Table~1 below.

\begin{center}
\begin{supertabular}{|c|c|c|c|}
\hline $r$ & $\dim_{\mathbb{R}} S^r_d(\Delta_1)$ & $\dim_{\mathbb{R}} S^r_d(\Delta_2)$ & $\dim_{\mathbb{R}} S^r_d(\Delta_3)$\\
\hline $0$ & $2d^2+2d+1$       & $2d^2+2d+1$     & $2d^2+2d+1$ \\
\hline $1$ & $2d^2-2d+4$   & $2d^2-2d+3$  & $2d^2-2d+3$ \\
\hline $2$ & $2d^2-6d+13$  & $2d^2-6d+11$ & $2d^2-6d+10$ \\
\hline $3$ & $2d^2-10d+28$  & $2d^2-10d+24$ & $2d^2-10d+23$ \\
\hline $4$ & $2d^2-14d+49$ & $2d^2-14d+43$ & $2d^2-14d+41$ \\
\hline
\end{supertabular}
\end{center}
\begin{center}
Table 1: Dimension of the spline spaces for $\Delta_1, \Delta_2, \Delta_3$
\end{center}

\noindent In \S2.3 we describe the code used to compute these dimensions. There are three things to highlight about Table~1.
\begin{itemize}
    \item First, for all three examples and all orders of smoothness $r$, the dimension of all the spline spaces is given by a polynomial function of $d$. This is true when $\Delta$ has a single interior vertex (Schumaker \cite{schumaker1}) but not in general. Schumaker's work is discussed in \S 2, which also includes the code used to produce the dimensions. \vskip .05in
\item Second, for all three examples, if we fix the order of smoothness $r$ then the polynomials $ad^2+bd+c$ computing the dimension of $S^r_d(\Delta_i)$ have the same coefficients $a$ and $b$. Schumaker gave a formula for these two coefficients (see Equation~\ref{bound} of \S 2) and for the coefficient $c$. As the reader might guess from this example, more symmetry in the triangulation leads to a larger space of splines.\vskip .05in
\item Third, note that $\dim_{\RR}S^0_d(\Delta_i)$ is the same for all three triangulations. It turns out that this is not an accident: the dimension of $S^0_d(\Delta_i)$ is determined solely by the combinatorics of the triangulation.
\end{itemize}
\end{example0}

\begin{example0}\label{ex:second}
The next example relates to a famous conjecture of Strang \cite{strang} and Billera's beautiful answer in \cite{bTAMS}, which we discuss in detail in the next section. Consider the two triangulations below.  As in Example~\ref{ex:first}, the combinatorics of the two examples are the same. Moreover, in this example, the local geometry at each of the three interior vertices is also the same, in the sense that there are four edges incident to each interior vertex, and at each vertex there are four distinct slopes. 
\begin{figure}[ht]
\vskip -.2in
\hskip -.6in\includegraphics[width=6in]{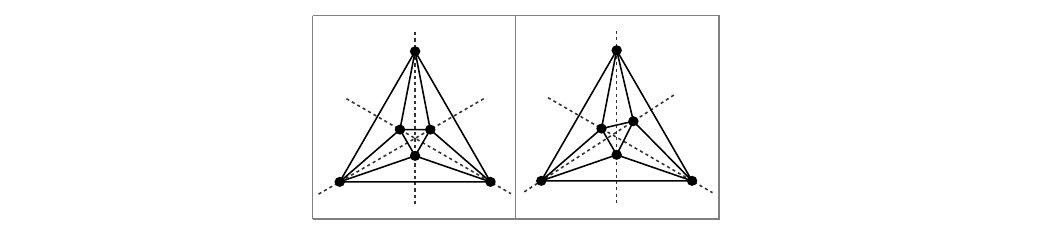}
\caption{Triangulations $\Delta_4$ and $\Delta_5$ with same combinatorics and local geometry.  (Dotted lines connect each boundary vertex with the interior vertex furthest from it but are not part of the triangulation.) } \label{figure: Morgan-Scott}
\end{figure}

\noindent However, there is subtle global geometry associated with these two triangulations: the dotted lines show us that $\Delta_4$ has more symmetry than $\Delta_5$. More precisely, these three dotted lines are incident at a point in $\Delta_4$ but not in $\Delta_5$. 
A computation shows that for all $d \ge 3$, the dimensions of the spline spaces satisfy~\[
\dim S^1_d(\Delta_4)=\frac{7}{2}d^2-\frac{15}{2}d+7 = \dim S^1_d(\Delta_5).
\]

\noindent In \cite{ms}, Morgan and Scott found that the splines behave specially when $d = 2$:
\[
\begin{array}{ccc}
\dim S^1_2(\Delta_4) & =& 7, \\
& & \\
\dim S^1_2(\Delta_5) & = & 6.
\end{array}
\]
In particular, while the dimension of the spline space is given by a polynomial depending only on combinatorics and local geometry when $d$ is sufficiently large, interesting things happen for small $d$. In \S 2, we will see that this is related to {\em homology}. For more analysis of this example, see Diener's work \cite{diener}.
\end{example0}

\subsection{Dualizing the definition of splines} 
Up to now, we've asserted without evidence that splines provide a bridge between the seemingly vast gulf between numerical analysis and equivariant cohomology. In this section, we describe a dual construction of splines.  From one point of view, this dual construction can be treated as the definition of equivariant cohomology. We sketch another point of view in \S 1.4, with more details in \S 3. 

The dual is constructed in two steps:
\begin{enumerate}
    \item {\bf Combinatorially construct a graph:} Suppose $\Delta$ is a polyhedral subdivision.  Construct the graph $G_{\Delta}^*$ that is dual in a graph-theoretic sense: create a vertex $v_i$ in $G_{\Delta}^*$ for each top-dimensional face $F_i$ in $\Delta$  and create an edge $v_iv_j$ between two vertices if and only if the corresponding faces $F_i \cap F_j$ intersect in a face of dimension one less than either $F_i$ or $F_j$. For instance if $\Delta$ is a triangulation in the plane, then for each triangle $T_i$ in $\Delta$ construct a vertex $v_i$ in the vertex set $V$ and make the rule that $v_iv_j$ is an (undirected) edge in $E$ if and only if $T_i \cap T_j$ is an edge of both corresponding triangles.
    \item {\bf Algebraically determine labels on the graph's edges:} Suppose that $F_i, F_j$ are two faces in $\Delta$ whose intersection $F_i \cap F_j$ has codimension one, namely the intersection is defined by an equation of the form
    \[\sum_\ell a_\ell x_\ell + c_{ij} = 0.
    \]
    Then label the edge $v_iv_j$ in $G_{\Delta}^*$ by 
    \[\ell(v_iv_j) = \sum_{\ell} a_{\ell} x_{\ell} +c_{ij}.\]
    We call the function $\ell: E \longrightarrow \mathbb{R}[x_1,\ldots,x_n]$ an  {\textit{edge-labeling}} of the graph $G_{\Delta}^*$.
\end{enumerate}

\begin{example0} \label{ex: dual graphs}
For example, the figure below 
shows the dual graph to the triangulation in the middle of Figure~\ref{figure: star triangulations}.  

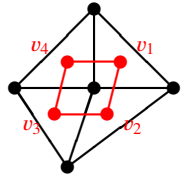
\begin{figure}[h]
\hskip 2in\begin{picture}(60,60)(0,0)
\thicklines
\put(30,60){\circle*{5}}
\put(0,30){\circle*{5}}
\put(30,30){\circle*{5}}
\put(60,30){\circle*{5}}
\put(20,0){\circle*{5}}
\put(20,0){\line(1,3){10}}
\put(20,0){\line(4,3){40}}
\put(20,0){\line(-2,3){20}}
\put(30,60){\line(-1,-1){30}}
\put(30,60){\line(1,-1){30}}
\put(0,30){\line(1,0){60}}
\put(30,30){\line(0,1){30}}

\put(20,40){\color{red} \circle*{5}}
\put(6,44){\color{red} {$v_4$}}
\put(40,40){\color{red} \circle*{5}}
\put(46,44){\color{red} {$v_1$}}
\put(15,20){\color{red} \circle*{5}}
\put(35,20){\color{red} \circle*{5}}
\put(3,14){\color{red} {$v_3$}}
\put(41,14){\color{red} {$v_2$}}

\put(20,40){\color{red} \line(1,0){20}}
\put(15,20){\color{red} \line(1,0){20}}
\put(15,20){\color{red} \line(1,4){5}}
\put(35,20){\color{red} \line(1,4){5}}
\end{picture}
\caption{The dual graph (in red) of a triangulation (in black)}
\end{figure} \label{figure: computing dual graph}

As a combinatorial object, the graph $G_{\Delta}^*$ only consists of a set of vertices and edges, not a particular embedding in the plane.  This means that the dual graph of each example in Figure~\ref{figure: star triangulations} is the same: a cycle on four vertices, like that shown on the left in Figure~\ref{figure: running example dual graphs with labelings} or (with additional labels) in Figure~\ref{figure: computing dual graph}.  In fact, Figure~\ref{figure: running example dual graphs with labelings} shows the graphs dual to the running examples in Figures~\ref{figure: star triangulations} and~\ref{figure: Morgan-Scott}.   

Where these examples differ is in their edge-labelings.  For instance, consider the star triangulations in Figure~\ref{figure: star triangulations}.  Assume without loss of generality that the point in the center is the origin $(0,0)$.  The four triangles in the leftmost example in Figure~\ref{figure: star triangulations} share either an edge on the line $y=0$ or an edge on the line $x=0$ so the edge-labels on the corresponding dual graph are
\[\ell(v_1v_2) = \ell(v_3v_4) = y, \hspace{0.5in} \textup{ and } \hspace{0.5in} \ell(v_1v_4) = \ell(v_2v_3) = x.\]
The second star differs from the first only in the bottom edge, which now lies on the line $-3x+y=0$.  Thus the edge-labels for the dual graph in this case are 
\[\ell(v_1v_2) = \ell(v_3v_4) = y, \qquad \ell(v_1v_4) = x  \qquad \textup{ and } \qquad \ell(v_2v_3) = -3x+y. \]
Finally, the rightmost triangulation is obtained from the middle one by stretching the triangles whose intersection defines the edge-label $\ell_{34}$.  That edge in the triangulation now lies on the line $3x+y=0$ and so the edge-labels are now
\[ \ell(v_1v_2)=y, \qquad  \ell(v_1v_4) = x, \qquad \textup{ and } \ell(v_2v_3) = -3x+y, \qquad \ell(v_3v_4) = 3x+y \]
\begin{figure}[h]
\hskip .8in\begin{picture}(60,60)(0,-10)
\thicklines 
\put(20,40){\circle*{5}}
\put(40,40){\circle*{5}}
\put(15,20){\circle*{5}}
\put(35,20){\circle*{5}}

\put(20,40){\line(1,0){20}}
\put(26,44){$\ell_{14}$}

\put(15,20){\line(1,0){20}}
\put(5,28){$\ell_{34}$}

\put(15,20){\line(1,4){5}}
\put(20,11){$\ell_{23}$}

\put(35,20){\line(1,4){5}}
\put(40,27){$\ell_{12}$}
\end{picture} \hspace{1in}
\begin{picture}(60,70)(-30,-40)
\thicklines 
\put(0,0){\circle*{5}}
\put(0,30){\circle*{5}}
\put(-30,-20){\circle*{5}}
\put(30,-20){\circle*{5}}
\put(30,15){\circle*{5}}
\put(-30,15){\circle*{5}}
\put(0,-35){\circle*{5}}

\put(0,0){\line(0,1){30}}
\put(0,0){\line(-3,-2){30}}
\put(0,0){\line(3,-2){30}}
\put(-11,10){$\ell_{01}$}
\put(-12,-15){$\ell_{05}$}
\put(12,-6){$\ell_{03}$}

\put(30,15){\line(0,-1){35}}
\put(30,15){\line(-2,1){30}}
\put(18,24){$\ell_{12}$}
\put(33,-5){$\ell_{23}$}

\put(-30,15){\line(0,-1){35}}
\put(-30,15){\line(2,1){30}}
\put(-24,26){$\ell_{16}$}
\put(-42,-2){$\ell_{56}$}

\put(0,-35){\line(-2,1){30}}
\put(0,-35){\line(2,1){30}}
\put(-23,-35){$\ell_{45}$}
\put(18,-35){$\ell_{34}$}
\end{picture}
\caption{Dual graphs for the triangulations in Figures~\ref{figure: star triangulations} and~\ref{figure: Morgan-Scott}: different geometries come from the choice of labels $\ell_{ij}=\ell(v_iv_j)$} \label{figure: running example dual graphs with labelings}
\end{figure}
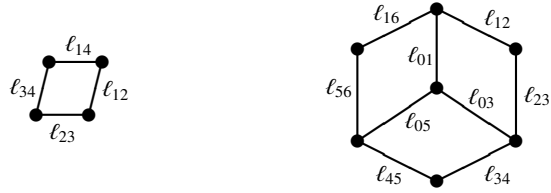
\end{example0}
\vskip -.1in
\noindent Let $\mathbb{R}[x_1\ldots, x_n]_{\leq d}$ denote the vector space of polynomials of degree at most $d$.
\begin{definition}\label{SplinesForGraphs}
    Let $G=(V,E,\ell)$ be a graph with vertex set $V$, edge set $E$ and an edge-labeling function $\ell:E\rightarrow \mathbb{R}[x_1, \ldots, x_n]$ for which the image of $\ell$  is contained in the subspace of (total) degree $1$ polynomials (not necessarily homogeneous). The space of splines $f_\bullet$ on $G$ of smoothness $r$ and degree at most $d$ is 
    \[
\begin{array}{l} S_d^r(G)= \\
\hspace{0.15in}\left\{f_{\bullet} = (f_v)\in\bigoplus\limits_{v \in V} \mathbb{R}[x_1, \ldots, x_n]_{\leq d}\Big |\ \begin{array}{c} \hbox{ for each }(vw)\in E,  \hbox{ the polynomial } \\ \ell(vw)^{r+1} \hbox{ divides } f_v-f_w \end{array} \right\}
 \end{array}   \]
 and also $S^r(G)=\bigcup_{d=0}^\infty S^r_d(G)$.
\end{definition}
It is clear that $S_d^r(G)$ is a vector space: if $f_\bullet$ and $g_\bullet$ satisfy the divisibility conditions imposed by the edge-labeling, then so does their sum $f_\bullet + g_\bullet$ and all scalar multiples $\lambda f_\bullet=(\lambda f_v)_{v\in V}$ for any $\lambda\in \mathbb{R}$. In fact, the divisibility conditions are preserved even when we multiply $f_\bullet$ and $g_\bullet$ coordinate wise, meaning that $S^r(G)$ has the additional structure of an algebra over the polynomial ring $\mathbb{R}[x_1, \ldots, x_n]$.

\begin{example0}
If $\Delta \subseteq \mathbb{R}$ consists of two intervals with intersection at the origin then its dual graph $G^*_\Delta$ consists of two vertices joined by a single edge.  The  intersection of the two faces of $\Delta$ is defined by the equation $x=0$ so the corresponding edge-label in $G^*_\Delta$ is $x$. 

Now suppose that $p_1$ and $p_2$ are polynomials on the two faces of $\Delta$.  They agree at $x=0$ if and only if the polynomial $p_1 - p_2$ has no constant term, or equivalently $x$ divides $p_1-p_2$. The derivatives of $p_1$ and $p_2$ agree at $x=0$ if and only if $p_1' - p_2'$ has no constant term, or equivalently $x^2$ divides $p_1 - p_2$, and similarly for all higher derivatives. So in this case, splines on $\Delta$ and splines on $G^*_\Delta$ both consist of the same pairs of polynomials.

Moreover, the claim is true regardless of the point of intersection of the two intervals in $\Delta$.  For instance,  if we changed the point of intersection of the two faces in $\Delta$ to $x=1$ then the corresponding edge-label in $G^*_\Delta$ would become $x-1$, but the rest of the argument would remain the same.
\end{example0}

\begin{example0}
 More generally, Billera proved that for most polydedral subdivisions $\Delta$ that arise in applications, the vector space $S_d^r(G^*_\Delta)$ coincides with the spline space $S_d^r(\Delta)$ defined before Equation~\ref{SplineModule}. Section~\ref{sec: homology} discusses this important result as well as its historical context in more detail. 
 \end{example0}
\subsubsection{Varieties and Torus actions} Now we switch gears.  Suppose $X$ is an algebraic variety---namely, the solution set of a system of polynomial equations---that is in addition equipped with an action of a torus $T$. Intuitively, the torus can be thought of as acting by rotations that can be related to flows along a vector field when $X$ is smooth.  While $T$ is often modeled by a product group of rotations in various planes (denoted $S^1 \times S^1 \times \cdots \times S^1$), it can be useful to instead represent $T$ as the \textit{algebraic torus}: the product group $\mathbb{C}^* \times \mathbb{C}^* \times \cdots \times \mathbb{C}^*$ with $n$ factors, so $T$ acts by pointwise multiplication on $\mathbb{C}^n$. On the one hand, the difference is small: if we write $\mathbb{C}^*$ using polar coordinates, we recover rotations by the circle $S^1$ just by ignoring dilations.  On the other hand, some people feel strongly about passing from $\mathbb{R}$ to $\mathbb{C}$, either for disciplinary reasons (symplectic geometers prefer $\mathbb{R}$ while algebraic geometers and algebraic topologists tend to prefer $\mathbb{C}$) or for mathematical reasons (it is easier to describe various algebraic calculations --- like roots of polynomials or eigenvalues of matrices --- in a uniform way over $\mathbb{C}$).  

We will assume that $X$ and its $T$-action are both ``well-behaved" in a technical sense not dissimilar to our assumptions on the partitioned domain $\Delta$ (see surveys like \cite{Ty2005} for more on these conditions). Then the set of $T$-fixed points and one-dimensional complex $T$-orbits define the vertices and edges of a graph $G_{X,T}$. Moreover, the torus action restricts to a rotation of the one-dimensional orbits, which we can encode as an edge-labeling $\ell$ of the graph $G_{X,T}$.  In this case, pioneering work of Goresky, Kottwitz and MacPherson \cite{GKM98} proved that the set of splines $S^0(G_{X,T,\ell})$ is isomorphic to the $T$-equivariant cohomology ring of $X$. We say more about this approach, which is often called \textit{GKM theory}, in \S 1.4.

\begin{example0} \label{example: P1 moment graph}
For example, consider the (complex) projective space $X=\mathbb{P}^1$.  Each element of $X$ is a line through the origin in $\mathbb{C}^2$, which we denote using any nonzero vector $[a,b]$ parametrizing the line. Now suppose $T=\mathbb{C}^* \times \mathbb{C}^*$ acts on $X$ via 
\[(z_1, z_2)\cdot \left[
\begin{array}{c}
a,
b 
\end{array}\right]=\left[
\begin{array}{c}
z_1 a,
z_2 b 
\end{array}\right].\]  
What it means for $[a,b]$ to be a $T$-fixed point is that 
\[(z_1,z_2)\cdot \left[\begin{array}{c}
a,
b 
\end{array}\right]=\left[\begin{array}{c}
a,
b 
\end{array}\right]\] 
for all $(z_1,z_2)\in T$.  This is true if and only if $\left[\begin{array}{c}
a,
b 
\end{array}\right]$ is either $\left[\begin{array}{c}
1,
0 
\end{array}\right]$ or $\left[\begin{array}{c}
0,
1 
\end{array}\right]$ since these are the only two vectors whose slope remains constant no matter which $z_1, z_2$ we choose. Since $X$ itself has one (complex) dimension, the only one-dimensional orbit consists of the set of points where both coordinates $a, b \neq 0$.  The closure of this one-dimensional orbit contains both of the fixed points.  Put together, this means the graph $G_{X,T}$ consists of one edge (corresponding to the single one-dimensional orbit) that connects the two fixed points. 

The edge label is given by the torus \textit{character}, which we can think of as a linearization of the action of $T$ on this orbit.  For instance, suppose we take the element of the one-dimensional orbit of $X$ whose slope is $a/b$.  Consider the slope of the image of this element under the action of $T$:
\[(z_1,z_2)\cdot \left[\begin{array}{c}
a,
b 
\end{array}\right]=\left[\begin{array}{c}
z_1 a,
z_2 b 
\end{array}\right].
\] 
In other words, the image has slope $(z_1 a)/(z_2 b) = (z_1/z_2) (a/b)$. If we restrict $T$ to this one-dimensional orbit, it acts as the one-dimensional torus parametrized by $z_1/z_2$ which we can write additively as the character $x_1 - x_2$.  (In fact, this character is uniquely defined up to sign.)

In short, with $X = \mathbb{P}^1$ and $T = \mathbb{C}^* \times \mathbb{C}^*$ as above, the graph $G_{X,T}$ has two fixed points $v_1 = [1,0]^T$ and $v_2 = [0,1]^T$ and one edge $v_1v_2$.  Moreover the edge-labeling is defined by $\ell(v_1v_2)=x_1-x_2$.
\end{example0}

\begin{example0}
The previous example extends without much difficulty to the case when $X = \mathbb{P}^n$ is general projective space. The elements of $\mathbb{P}^{n}$ consist of lines through the origin in $\mathbb{C}^{n+1}$, which we'll denote (non-uniquely) by any choice of vector $[a_1,\ldots,a_{n+1}]$ spanning a particular line.  Let $T = (\mathbb{C}^*)^{n+1}$ act on $X$ by
\[(z_1, z_2,\ldots, z_{n+1}) \cdot [a_1, a_2, \ldots, a_n,a_{n+1}] = [z_1a_1,z_2a_2,\ldots,z_{n+1}a_{n+1}].
\]
The analysis of $\mathbb{P}^1$ in Example~\ref{example: P1 moment graph} shows that if $i, j$ is any pair of entries that are both nonzero then some elements of $T$ change the ratio $a_i/a_j$.  This means the only elements of $X$ that are fixed by \textit{all} choices of $(z_1, z_2, \ldots, z_{n+1})$ are the coordinate lines, namely those for which exactly one of the $a_i$ is nonzero.  Thus $G_{X,T}$ has $n+1$ vertices $v_1, \ldots, v_{n+1}$ corresponding to the $n+1$ coordinate axes of $\mathbb{C}^{n+1}$.  

Similarly, the one-dimensional (complex) orbits of $T$ are characterized by the choice of exactly \textit{two} nonzero entries $i, j$ and the closure of each one-dimensional $T$-orbit contains the $T$-fixed points corresponding to the $i^{th}$ and $j^{th}$ coordinate axes.  Thus $G_{X,T}$ has edges $v_iv_j$ for each pair $i \neq j$ of distinct vertices --- in other words, the graph $G_{X,T}$ is a complete graph on $n+1$ vertices.
\pagebreak

Finally, suppose $i, j$ are the two nonzero entries in a one-dimensional orbit.  The action of $(z_1,z_2, \ldots, z_{n+1})$ on an element $[0,\ldots, 0, a_i, 0,\ldots, 0,a_j,0,\ldots,0]$ of the one-dimensional orbit changes the ratio $a_i/a_j$ to the ratio $(z_i a_i)/(z_j a_j)$.  In other words, the torus action on this one-dimensional orbit scales the ratio $a_i/a_j$ by $z_i/z_j$.  Linearizing gives the edge-label $\ell(v_iv_j) = x_i - x_j$, which is the character of the torus action on this one-dimensional orbit. There is an interesting connection between this example and the triangulation called the Alfeld split; for more, see \cite{s2, s3}.
\end{example0}
More details and examples can be found in surveys on GKM theory, such as \cite{Ty2005}.

It is often convenient to denote a spline on an edge-labeled graph by labeling the vertices of the graph with their corresponding polynomials. For example, if $G$ is the first graph in Fig. \ref{figure: running example dual graphs with labelings} with labels coming from the second triangulation of Example \ref{ex:first}, an element of $S_1^0(G)$ is given in Figure~\ref{figure: spline element} below. 

\vskip -.2in
\begin{figure}[h]
\hskip 1.5in \begin{picture}(60,100)(-20,0)
\thicklines 
\put(33,83){$2x$}
\put(78,83){$-x$}
\put(20,32){$2x+y$}
\put(65,32){$-x+2y$}
\put(40,80){\line(1,0){40}}

\put(30,40){\line(1,0){40}}

\put(30,40){\line(1,4){10}}

\put(70,40){\line(1,4){10}}
\put(110,55){$\in S^0_1(G)$}
\end{picture} 
\vskip -.4in
\caption{An element of the spline for the graph coming from the second triangulation of Example \ref{ex:first}.}
\label{figure: spline element}
\vskip -.15in
\end{figure}
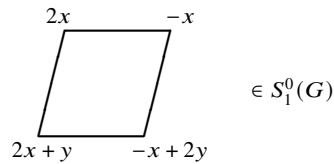
\subsubsection{Finding bases for splines on graphs}
One core problem in equivariant cohomology is finding bases for splines on a graph, which in this context means ``giving a closed formula for all the polynomials of each basis spline."  The most common kind of basis in geometric and topological contexts is an analogue of upper-triangular bases.  From one perspective, these bases arise from choosing a north pole for the variety, using the torus action to flow points on the variety towards the north pole, and then observing in what subvarieties the points get stuck.  This is the idea of Morse theory \cite{milnor} in topology, adapted by Bialynicki-Birula to complex algebraic varieties in \cite{BB1, BB2}.

Motivated by this philosophy, we can construct the same bases purely algebraically via a sort of Gaussian elimination, as done in, e.g., \cite{GTolman, KT, Ty08b}.  Rather than formally describing this process, we give two illustrative examples.

\begin{example0}\label{ex: basis of spline as a module}Consider the edge-labeled graph $G$ in Figure~\ref{figure : graph of Delta 1}, corresponding to the triangulation $\Delta_1$ in Figure~\ref{figure: star triangulations}.
\vskip -.3in
\begin{figure}[ht]
\hskip 1.4in \begin{picture}(60,80)(-20,20)
\thicklines 
\put(33,82){$v_4$}
\put(81,82){$v_1$}
\put(21,36){$v_3$}
\put(72,36){$v_2$}

\put(40,80){\line(1,0){40}}
\put(56,84){$x$}

\put(30,40){\line(1,0){40}}
\put(26,60){$y$}

\put(30,40){\line(1,4){10}}
\put(49,33){$x$}

\put(70,40){\line(1,4){10}}
\put(80,60){$y$}
\end{picture} 
\vskip -.2in
\caption{Graph corresponding to the triangulation $\Delta_1$ in Figure~\ref{3planarTriangulations}} \label{figure : graph of Delta 1}
\vskip -.2in
\end{figure}
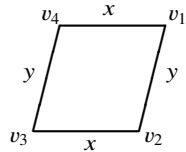

\noindent We represent an element of the corresponding spline space $S^r(G)$ as \[f_\bullet=(f_{v_1}, f_{v_2}, f_{v_3}, f_{v_4}).\] The conditions given by the edge labels are
\[
\begin{array}{c}
y^{r+1} \mid f_{v_1}-f_{v_2}\\
x^{r+1} \mid f_{v_2}-f_{v_3}\\
y^{r+1} \mid f_{v_3}-f_{v_4}\\
x^{r+1} \mid f_{v_1}-f_{v_4}.
\end{array}
\]
These conditions are equivalent to 
\begin{equation}\label{eqn : spline basis element example}
\begin{array}{ccc}
f_{v_2}&=&f_{v_1}+a y^{r+1}\\
f_{v_3}&=&f_{v_2}+ bx^{r+1}\\ 
f_{v_4}&=& f_{v_1}+cx^{r+1} \\
 & = &f_{v_3}+dy^{r+1}   
 \end{array}
\end{equation}
 for polynomials $a,b,c,d\in\mathbb{R}[x,y]$. 
 By substituting the expression for $f_{v_2}$ into the equation for $f_{v_3}$ and similarly for $f_{v_3}$ into $f_{v_4}$, we obtain 
\[
\begin{array}{ll}
f_{v_3} =f_{v_1}+a y^{r+1}+ bx^{r+1}  \\
f_{v_4} = f_{v_1}+a y^{r+1}+ bx^{r+1} + dy^{r+1}
\end{array}
\]
Since $f_{v_4}-f_{v_1}$ also equals $bx^{r+1}$, we deduce that the polynomials $a,b,c,d$ solve the system of Equations \eqref{eqn : spline basis element example} if and only if 
\[
y^{r+1}(a+d) + x^{r+1}b = cx^{r+1}.\]
Since $y^{r+1}$ and $x^{r+1}$ are relatively prime, we can rephrase the above condition as
\[
d=-a + x^{r+1}u \quad \textup{     and  }
\quad c=b+uy^{r+1}
\]
for some polynomial $u\in \mathbb{R}[x,y]$. In other words, any triple of polynomials $a, b, u \in \mathbb{R}[x,y]$ gives a unique spline $f_\bullet \in S^r(G)$ and vice versa. 

Figures \ref{figure: spline basis Gaussian elimination 1}, \ref{figure: spline basis Gaussian elimination 2}, and \ref{figure: spline basis Gaussian elimination 3} graphically illustrate what we just discussed. By looking at the pictures, it should be clear why we call the process a sort of Gaussian elimination.
\begin{figure}[h]
\scalebox{1.25}{
\begin{picture}(100,100)(20,0)
\thicklines 
\put(30,83){$f_{v_4}$}
\put(77,83){$f_{v_1}$}
\put(19,30){$f_{v_3}$}
\put(67,30){$f_{v_2}$}
\put(37,80){\line(1,0){40}}

\put(27,40){\line(1,0){40}}

\put(27,40){\line(1,4){10}}

\put(67,40){\line(1,4){10}}
\put(97,55){$=$}

\put(120,83){$f_{v_1}$}
\put(167,83){$f_{v_1}$}
\put(109,30){$f_{v_1}$}
\put(157,30){$f_{v_1}$}
\put(127,80){\line(1,0){40}}

\put(117,40){\line(1,0){40}}

\put(117,40){\line(1,4){10}}

\put(157,40){\line(1,4){10}}
\put(185,55){$+$}
\put(200,94){$f_{v_4}-f_{v_1}$}
\put(195,84){\tiny{ $=bx^{r+1}+ux^{r+1}y^{r+1}$}
}
\put(260,85){$0$}
\put(190,30){$f_{v_3}-f_{v_1}$}
\put(190,18){\tiny{ $=ay^{r+1}+bx^{r+1}$}
}
\put(245,30){
$f_{v_2}-f_{v_1}
$
}
\put(245,18){\tiny{ $=ay^{r+1}$}}
\put(220,80){\line(1,0){40}}

\put(210,40){\line(1,0){40}}

\put(210,40){\line(1,4){10}}

\put(250,40){\line(1,4){10}}

\end{picture} 
}
\vskip -.2in
\caption{First step of \emph{Gaussian elimination} like process, with computed values for $f_{v_i}$.}
\label{figure: spline basis Gaussian elimination 1}
\end{figure}

\end{example0}

\begin{figure}[ht]
\vskip -.1in
\scalebox{1.25}{
\begin{picture}(100,100)(0,0)
\thicklines 
\put(14,85){\tiny{$bx^{r+1}+ux^{r+1}y^{r+1}$}}
\put(72,83){$0$}
\put(0,30){$ay^{r+1}+bx^{r+1}$}
\put(62,30){$ay^{r+1}$}
\put(33,80){\line(1,0){40}}

\put(23,40){\line(1,0){40}}

\put(23,40){\line(1,4){10}}

\put(63,40){\line(1,4){10}}
\put(93,55){$=$}

\put(120,83){$0$}
\put(167,83){$0$}
\put(109,30){$ay^{r+1}$}
\put(157,30){$ay^{r+1}$}
\put(127,80){\line(1,0){40}}

\put(117,40){\line(1,0){40}}

\put(117,40){\line(1,4){10}}

\put(157,40){\line(1,4){10}}
\put(185,55){$+$}
\put(200,85){\tiny{$bx^{r+1}+ux^{r+1}y^{r+1}$}}
\put(260,83){$0$}
\put(250,30){$0$}
\put(205,30){$bx^{r+1}$}
\put(220,80){\line(1,0){40}}

\put(210,40){\line(1,0){40}}

\put(210,40){\line(1,4){10}}

\put(250,40){\line(1,4){10}}

\end{picture} }
\vskip -.2in
\caption{Second step of \emph{Gaussian elimination} like process.}
\label{figure: spline basis Gaussian elimination 2}
\end{figure}
\begin{figure}[ht]
\scalebox{1.25}{
\begin{picture}(100,100)(0,0)
\thicklines 
\put(14,85){\tiny{$bx^{r+1}+ux^{r+1}y^{r+1}$}}
\put(70,83){$0$}
\put(11,30){$bx^{r+1}$}
\put(60,30){$0$}
\put(33,80){\line(1,0){40}}

\put(23,40){\line(1,0){40}}

\put(23,40){\line(1,4){10}}

\put(63,40){\line(1,4){10}}
\put(94,55){$=$}

\put(120,83){$bx^{r+1}$}
\put(167,83){$0$}
\put(109,30){$bx^{r+1}$}
\put(157,30){$0$}
\put(127,80){\line(1,0){40}}

\put(117,40){\line(1,0){40}}

\put(117,40){\line(1,4){10}}

\put(157,40){\line(1,4){10}}
\put(185,55){$+$}
\put(213,83){$ux^{r+1}y^{r+1}$}
\put(260,83){$0$}
\put(200,30){$0$}
\put(250,30){$0$}
\put(220,80){\line(1,0){40}}

\put(210,40){\line(1,0){40}}

\put(210,40){\line(1,4){10}}

\put(250,40){\line(1,4){10}}

\end{picture} }
\caption{Third (and final) step of \emph{Gaussian elimination} like process.}
\label{figure: spline basis Gaussian elimination 3}
\end{figure}

\noindent Now let
\[
g_\bullet^{(1)}=(1,1,1,1), \quad g_\bullet^{(2)}=(0,y^{r+1},y^{r+1},0),
\]
\[g_\bullet^{(3)}=(0,0,x^{r+1},x^{r+1}), \quad g_\bullet^{(4)}=(0,0,0,y^{r+1}x^{r+1}).
\]
Our previous calculations showed that $f_\bullet\in S^r(G)$ if and only if there exist $h_1, h_2, h_3, h_4\in\mathbb{R}[x,y]$ such that 
\[
f_\bullet= h_1 g_\bullet^{(1)}+ h_2 g_\bullet^{(2)}+ h_3 g_\bullet^{(3)} + h_4 g_\bullet^{(4)}.\] Moreover, it is easy to show that if the $h_i$ exist then they are unique. In other words, the $g_\bullet^{(i)}$ form a basis for $S^r(G)$ as a module over $\mathbb{R}[x,y]$.

The quest for bases with properties useful for illuminating intrinsic structure has been a central theme in GKM theory, see for example \cite{Go14, HTy17, LP20}.

\begin{example0} \label{example: degree-two splines}
We can compute bases for the other examples in Figure~\ref{figure: star triangulations} using exactly the same strategy, but changing the edge-labeling will change the results.  For instance, consider the edge-labeling in Figure~10 below, which corresponds to the triangulation $\Delta_2$ in Figure~\ref{3planarTriangulations}.
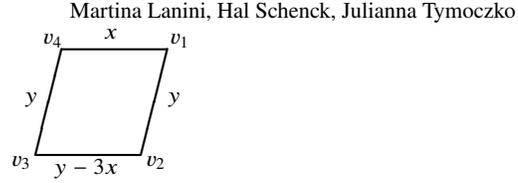
\begin{figure}[h]\label{figg10}
\hskip 1.4in \begin{picture}(60,50)(-20,30)
\thicklines 
\put(33,82){$v_4$}
\put(81,82){$v_1$}
\put(21,36){$v_3$}
\put(72,36){$v_2$}

\put(40,80){\line(1,0){40}}
\put(56,84){$x$}

\put(30,40){\line(1,0){40}}
\put(26,60){$y$}

\put(30,40){\line(1,4){10}}
\put(37,33){$y-3x$}

\put(70,40){\line(1,4){10}}
\put(80,60){$y$}
\end{picture} 
\label{figure : graph of Delta 2}
\caption{Graph corresponding to the triangulation $\Delta_2$ in Figure~\ref{3planarTriangulations}}
\end{figure} 

\noindent For computational convenience, we do this calculation in the case $r=1$. Suppose $f_{\bullet}$ is an arbitrary element of $S^1(G)$ and as before consider 
\[f_{\bullet} - (f_{v_1}, f_{v_1}, f_{v_1}, f_{v_1}) = (0,f_{v_2}^{(1)}, f_{v_3}^{(1)}, f_{v_4}^{(1)})\] as shown in Figure~\ref{figure: second spline basis, Gaussian elimination 1} below.
\begin{figure}[h]
\begin{picture}(100,100)(20,0)
\thicklines 
\put(30,83){$f_{v_4}$}
\put(77,83){$f_{v_1}$}
\put(19,30){$f_{v_3}$}
\put(67,30){$f_{v_2}$}
\put(37,80){\line(1,0){40}}

\put(27,40){\line(1,0){40}}

\put(27,40){\line(1,4){10}}

\put(67,40){\line(1,4){10}}
\put(97,55){$=$}

\put(120,83){$f_{v_1}$}
\put(167,83){$f_{v_1}$}
\put(109,30){$f_{v_1}$}
\put(157,30){$f_{v_1}$}
\put(127,80){\line(1,0){40}}

\put(117,40){\line(1,0){40}}

\put(117,40){\line(1,4){10}}

\put(157,40){\line(1,4){10}}
\put(185,55){$+$}
\put(213,84){$f^{(1)}_{v_4}$}
\put(260,83){$0$}
\put(203,30){$f^{(1)}_{v_3}$}
\put(250,30){$f^{(1)}_{v_2}$}
\put(220,80){\line(1,0){40}}

\put(210,40){\line(1,0){40}}

\put(210,40){\line(1,4){10}}

\put(250,40){\line(1,4){10}}

\end{picture} 
\caption{After one step in the \emph{Gaussian elimination} like process.}
\label{figure: second spline basis, Gaussian elimination 1}
\end{figure}
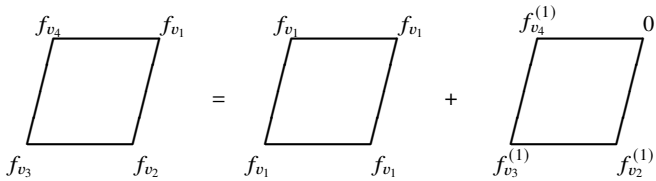
The polynomials $f^{(1)}_{v_i}$ satisfy similar conditions as before, except this time we have three distinct edge-labels:
\begin{equation}\label{eqn : second spline basis element example}
 f^{(1)}_{v_2}=a y^{2}, \quad
f^{(1)}_{v_3}=f^{(1)}_{v_2}+ b(y-3x)^{2}, \quad f^{(1)}_{v_4}= cx^2 = f^{(1)}_{v_3}+dy^{2}.   
\end{equation}
Note that we can subtract a polynomial multiple of $g^{(2)}=(0,y^2,y^2,0)$ to eliminate the entry at $v_2$. This gives us a spline of the form $(0,0,f_{v_3}^{(2)}, f_{v_4}^{(2)})$. Recall that the entries must satisfy 
\[
f_{v_3}^{(2)}=b(y-3x)^2 \quad \textup{ and } \quad f_{v_4}^{(2)}=cx^2=f_{v_3}^{(2)}+dy^2.
\]
Again, we can substitute $f^{(2)}_{v_3}$ into $f^{(2)}_{v_4}$ to get the equation
\[
cx^2 = b(y-3x)^2+dy^2.
\]
This time, note that $(y-3x)^2$ is not an $\mathbb{R}$-linear combination of $x^2$ and $y^2$ so no scalars $b,c, d \in \mathbb{R}$ solve this equation except $b=c=d=0$.  However, using a judicious choice of scalars $r_1, r_2, s_1, s_2 \in \mathbb{R}$, we can express any degree-three polynomial (uniquely) as a linear combination $(r_1x+r_2y)x^2 + (s_1x+s_2y)y^2$.  For instance, we have
\[
x(y-3x)^2 = xy^2 + (9x-6y)x^2, \quad \textup{ and } \quad y(y-3x)^2 = (-6x+y)y^2 + (9y)x^2.
\]
This means that amongst the splines on this graph that vanish simultaneously at $v_1$ and $v_2$, there are two $\mathbb{C}$-linearly independent splines of minimal degree, namely
\[
g^{(3a)} = (0,0,x(y-3x)^2,x^2(9x-6y)) \quad \hbox{
and }
\quad g^{(3b)} = (0,0,y(y-3x)^2,9yx^2). \]
This contrasts with the previous case, where we only had one basis vector $g^{(3)}$ for the vertex $v_3$.  Writing $f^{(2)}_{v_3}$ as a linear combination $a_1x(y-3x)^2+a_2y(y-3x)^2$ for some $a_1, a_2 \in \mathbb{R}[x,y]$, we use $g^{(3a)}$ and $g^{(3b)}$ to eliminate the entry $f^{(2)}_{v_3}$. What remains is a spline that is zero at all vertices except $v_4$, and this spline is clearly some multiple of \[g^{(4)} = (0,0,0,x^2y^2).\]

On the one hand, this process constructs five generators $g^{(1)}, g^{(2)}, g^{(3a)}, g^{(3b)}, g^{(4)}$ for splines on this edge-labeled graph.  That's one more generator than in the previous example!  On the other hand, the five splines are not linearly independent over \textit{polynomials}:
\[\begin{array}{ll}yg^{(3a)}-xg^{(3b)} &= (0,0,xy(y-3x)^2, x^2y(9x-6y))- (0,0,xy(y-3x)^2, 9x^3y)\\ &= (0,0,0,-6x^2y^2) \\ &= -6g^{(4)} \end{array}\]
So actually the first four splines generate $g^{(4)}$ over polynomials.  \end{example0}

The rightmost example $\Delta_3$ in Figure~\ref{3planarTriangulations} appears like it could be really different yet again, but in fact it is not. The reason is that if $a_1, a_2, a_3$ are any three \textit{distinct} scalars then $(x-a_1y)^2, (x-a_2y)^2, (x-a_3y)^2$ form a basis for the vector space of homogeneous degree-two polynomials, and if $a_1, a_2$ are any two distinct scalars then $x(x-a_1y)^2, y(x-a_1y)^2, x(x-a_2y)^2, y(x-a_2y)^2$ form a basis for the vector space of homogeneous degree-three polynomials.  This means the argument in Example~\ref{example: degree-two splines} holds whenever there are at least three distinct edge-labels in the graph --- and moreover, the argument holds whether the vector space is over $\mathbb{R}$ or $\mathbb{C}$.

Once we find a basis of $S^r(G)$ as an $\mathbb{R}[x_1, \ldots,x_n]$-module, determining the dimension of $S^r_d(G)$ for any $d$ becomes a combinatorial problem.

\begin{example0} Let $G$ be the edge-labeled graph $G$ in Figure~\ref{figure : graph of Delta 1}, corresponding to the triangulation $\Delta_1$ in Figure~\ref{figure: star triangulations}.
In Example \ref{ex: basis of spline as a module}, we found that for any $d\geq 0$ the dimension of $S^r_d(G)$
coincides with the dimension of the following $\mathbb{R}$-vector space:
\[
\left
\{
(h_1, h_2,h_3, h_4)\mid h_1\!\in \!\mathbb{R}[x,y]_{\leq d}, \ h_2,h_4 \!\in \!\mathbb{R}[x,y]_{\leq d-(r+1)}, \ h_3 \!\in\! \mathbb{R}[x,y]_{\leq d-2(r+1)}
\right\}
,\]
where $\dim(\mathbb{R}[x,y]_{\leq l})=0$ whenever $l<0$. 
To compute the dimension of $S^r_d(G)$ as an $\mathbb{R}$-vector space, we need only to incorporate the number of monomials of each degree.  This gives the following dimension formula:
\[
\dim(S_d^r(G))=
\left\{
\begin{array}{ll}
     {d+2 \choose 2}&\hbox{ if }d<r+1,  \vspace{0.5em} \\ 
    {d+2 \choose 2}+2{d-r+1\choose 2} & \hbox{ if }r+1\leq d<2(r+1), \vspace{0.5em} \\ 
     {d+2 \choose 2}+2{d-r+1\choose 2}+{d-2r\choose 2} & \hbox{ if }2(r+1)\leq d.
\end{array}
\right.
\]

\end{example0}

\subsection{Homology, Group Actions, Equivariant Cohomology}
To recap, \S 1.2 defined splines by gluing two functions together along a common face, subject to the compatibility conditions in Definition~\ref{SplineDef}; then \S 1.3 described a dual construction using labelled graphs and showed that it can be used to recover essentially the same information.  So why bother?  For two historical reasons, one algebraic reason, and one geometric reason.
\begin{enumerate}
    \item {\bf First historical reason:} In Section 2, we describe how Billera used this dual construction together with homological methods to prove Strang's conjecture.
    \vskip .1in
    \item {\bf Second historical reason:} Splines provide a model for equivariant cohomology of a large family of interesting spaces.  Some geometers define equivariant cohomology using piecewise polynomials on combinatorial subdivisions of a region in Euclidean space, especially algebraic geometers, and especially in the context of toric varieties, see for example \cite{brion, bVergne, payne, s2}.  However, the dual construction in this section is more common in symplectic geometry and algebraic topology.
    \vskip .1in
    \item {\bf Algebraic reason:}  From an algebraic perspective, the construction in this section is more natural --- or at least more generalizable, e.g. to arbitrary graphs $G$; with arbitrary coefficients (not just real or even complex polynomials but, say, finite fields); and with more complicated edge-labelings (not just by multiples of a single polynomial, but by more general ideals in the coefficients). Moreover, geometers and topologists care deeply about the ring structure of $S^r(G)$, which can be extended to a ring structure on $S^r_d(G)$ by using a quotient of the polynomial ring instead of the polynomial ring itself. 
\end{enumerate}
The last reason is the topic of \S 2 (on homology) and \S 3 (on group actions), which together describe two of the main algebraic tools used to study splines. This section provides a brief, intuitive introduction to the core ideas.

The concept of homology appeared more than a century ago, introduced by topologists interested in creating algebraic invariants capable of distinguishing between different topological spaces. The idea is to 
model a topological space $X$ by inductively gluing together combinatorial building blocks. Typically, the building blocks themselves are topologically uninteresting, but their intersections and unions can be quite complicated.  For example, imagine sections of fencing (topologically uninteresting, because they are just rectangles), which are connected together to form a corral (topologically interesting, because it divides space into two regions with an inside and an outside). 

The combinatorial construction defines a sequence of vector spaces---the {\em homology} of $X$---that are the desired algebraic invariants. Now, vector spaces have duals, and dualizing the construction above gives dual vector spaces. Perhaps surprisingly, the dual construction leads to richer algebra: the dual construction has the structure of a ring.  

It turns out that there are many variants of the construction defined above; the homology theory we use---simplicial homology---is described in detail in \S 2. Other variants involve adornments on either the space $X$ or the kinds of functions considered in the construction of homology.  For instance, one common situation in topology is when the space $X$ has some kind of symmetry, more precisely a \textit{group action}. This richer situation leads to  correspondingly more interesting  algebraic structure, the {\em equivariant cohomology}, which is the topic of \S 3. The punchline is that the ``dual" description of splines--via labelled graphs--matches perfectly into the setup of equivariant cohomology.

\section{Homology}\label{sec: homology}
As noted in the previous section, one problem that has attracted much attention is \textit{the dimension problem}: given a simplicial complex $\Delta \subset \RR^k$, determine the dimension of $\Srd$. Although it may seem hard to believe, this problem is open for $S^1_3(\Delta)$ when $\Delta \subseteq \RR^2$ is a triangulation in the plane.

We begin in \S 2.1 with the history of the dimension problem. In \S 2.2 we introduce {\em homology} formally and show how to use theoretical tools from homology to analyze questions about splines. In \S 2.3 we detour into computational algebra, demonstrating the {\tt AlgebraicSplines} package. Open problems and a discussion of current research directions appear in \S 4.1.

\subsection{Dimension of the spline space: history}
In 1974, Strang~\cite{strang} conjectured a formula for the
dimension of $S^1_d(\Delta)$ for $\Delta \subseteq \RR^2$. Schumaker proved a general lower bound for the dimension of spline spaces in \cite{schumaker1}, showing that for {\em any} planar triangulation $\Delta$ and all $d$ and $r$, the dimension of $\Srd$ is bounded below by the expression 
\begin{equation}\label{bound} 
L^r_d(\Delta) =  {d+2 \choose 2} + {d-r+1 \choose 2}f^0_1 - \left( {d+2 \choose 2} - {r+2 \choose 2} \right) f^0_0 +\gamma,
\end{equation}
\noindent with
\begin{equation}\label{boundLocal} 
\gamma = \sum_{v_i \in \Delta_0^0}\gamma_i \quad\mbox{ and }\quad\gamma_i = \sum_j \max\{(r+1+j(1-n(v_i))), 0 \}.
\end{equation}
In Equation~\ref{boundLocal}, $\Delta_0^0$ is the set of interior vertices and $n(v_i)$ is the number of distinct slopes at each interior vertex $v_i \in \Delta_0^0$.  This is a landmark result in the field, which Schumaker followed up with a complementary upper bound that differs from the lower bound only in the way the $\gamma_i$ terms are defined~\cite{schumaker2}.

\begin{example0}\label{dimensionFormulas} The triangulations of Example~\ref{ex:first} each have one interior vertex, and the triangulations differ only in the number of distinct slopes at that vertex, so Equation~\eqref{boundLocal} simplifies considerably.  In this case, the lower bound of Equation~\ref{bound} is equal to the results computed in Table 1. More generally, Schumaker's work shows that Equation~\ref{bound} gives the correct dimension whenever $\Delta$ is the star of a vertex. The terms $\gamma_i$ account for the local geometry. 
\end{example0}

\noindent Using Bernstein-B\`ezier methods, Alfeld and Schumaker~\cite{as2} proved that 
\[
L^r_d(\Delta) = \dim_{\RR}\Srd \mbox{ for }d\ge 4r+1.
\]
This was later extended to $d\geq 3r+2$ by Hong in~\cite{h}. 

How does this relate to Strang's conjecture? As Billera notes in \cite{bTAMS}, Strang was aware of the influence of local geometry on the dimension of the spline space. Strang conjectured that for a {\em generic} triangulation $\Delta$---a triangulation for which a miniscule perturbation of the vertices leaves the dimension unchanged---the dimension for the spline space is precisely $L^1_d(\Delta)$. Billera proved Strang's conjecture in \cite{bTAMS}.  In \S 2.2 we describe the homological approach that Billera used. 

\subsection{Homology and splines}
Solving the system of equations that encode the compatibility conditions in Definition~\ref{SplineDef} leads us inexorably to the homological setup. 

The first point to make is that {\em smoothness is a local condition}. In particular, suppose $\Delta \subseteq \RR^k$ is a cell complex and $\sigma_1$ and $\sigma_2$ are two $k$-cells sharing a common $k-1$ 
face $\tau=\sigma_1 \cap \sigma_2$.  Let $l_\tau$ be a (nonzero) polynomial of minimal degree vanishing on $\tau$. For example, if $\Delta$ is a simplicial or polyhedral subdivision then $l_\tau$ is a linear form, while for the semialgebraic splines of \S 4.1.7 the degree of $l_\tau$ can be higher. 

With this notation, it is a small exercise to show that if $F$ is a function on $\Delta$ then the restrictions $F|_{\sigma_1}$ and $F|_{\sigma_2}$ meet with smoothness $r$ across the shared face $\tau$ if and only if 
\begin{equation}\label{compatibility}
l_\tau^{r+1} \Big| \mbox{ }(F|_{\sigma_1}-F|_{\sigma_2}), \mbox{ where } l_\tau \mbox{ is a nonzero linear form vanishing on } \tau.
\end{equation}

\begin{example0}\label{ex:splineConditionsEx1} Consider the simplex $\Delta_1$ in Figure~1 of Example~\ref{ex:first}. Starting with the triangle in the first quadrant and moving clockwise, let $f_i=F|_{\sigma_i}$.  To obtain a spline of smoothness order $r$, we must solve the equations
\begin{center}
\vskip -.2in\begin{equation}\label{firstE}
\begin{array}{ccc}
a_1 y^{r+1}&=&f_1-f_2\\
a_2 x^{r+1}&=&f_2-f_3\\
a_3 y^{r+1}&=&f_3-f_4\\
a_4 x^{r+1}&=&f_4-f_1
\end{array}
\end{equation}
\end{center}
\end{example0}
Summing these equations yields 
\[a_1y^{r+1} + a_2x^{r+1} + a_3y^{r+1}+a_4x^{r+1} = 0\]
which generalizes to the other complexes in Example~\ref{ex:first} as $\sum_{i=1}^4a_i \ell_{i,i+1}^{r+1} = 0$ if we use the appropriate edge-labels $\ell_{i,i+1}$.  These are examples of a {\em syzygy}: for a set of polynomials $\{g_1, \ldots, g_m\}$, a syzygy is simply a polynomial relation: 
\[
\sum\limits_{i=1}^m a_ig_i = 0, \mbox{ where the }a_i \mbox{ are also polynomials}.
\]
Although far from obvious, the term $\gamma_i$ in Equation~\ref{boundLocal} encodes the dimension of the space of syzygies.  We will use this level of abstraction shortly. First we note that solving System~\ref{firstE} can also be viewed as computing the kernel of a matrix
\[
\phi =\left[ \begin{array}{rrrrcccc}
1&-1&0&0 &y^{r+1}&0&0&0\\
0&1&-1&0 & 0 &x^{r+1}&0&0\\
0&0 &1&-1 & 0& 0 &y^{r+1}&0\\
-1&0&0 &1 &0&0& 0& x^{r+1}
\end{array} \right] 
\]
It turns out that the $4 \times 4$ submatrix of $\pm 1$ on the left is secretly related to {\em homology}. To understand this, we need a primer on simplicial complexes and simplicial homology, for which we take a quick algebraic detour.
\subsubsection{A useful bookkeeping device: graded algebra}\label{GA}
The most natural way to express a spline on $\Delta \subseteq \RR^k$ is as a vector with one entry for each maximal face of $\Delta$, where each entry is a polynomial that satisfies certain smoothness conditions across $k-1$ dimensional faces of $\Delta$. As we saw in Equation~\ref{isModule}, the spline space $\Srd$ has the structure of a module over the polynomial ring $\RR[x_1,\ldots, x_k]$. 
It will be helpful for our bookkeeping to introduce {\em graded} rings and modules. The added notational complexity can be confusing but this extra work will pay off later.  
\begin{definition}\label{gradedObject}
For any integer $i \geq 0$, a \textit{homogeneous} polynomial of degree $i$ is a sum of monomials that all have degree $i$. Let $R_i$ denote the set of homogeneous polynomials of degree $i$.  Every polynomial ring $R$ is a direct sum
\[
R = \bigoplus_{i \in \ZZ}R_i
\]
and if $r_i \in R_i$ and $r_j \in R_j$ then $r_i \cdot r_j \in R_{i+j}$. These conditions make $R$ a ring \textit{graded by $\mathbb{Z}$}.  Similarly, a graded $R$-module $M =  \bigoplus_{i \in \ZZ} M_i$ is defined so if $r_i \in R_i$ and $m \in M_j$ then $rm \in M_{i+j}$. 

In our setting, $R_0=\RR$ is a field and therefore each $R_i$ is also
a vector space.  If $M$ is a graded module, then each $M_i$ is also an $R_0$ vector space.
\end{definition}

\begin{example0}\label{exm:shift}
Let $R=\RR[x,y]$. Then $x^2+3xy$ is homogeneous but $x^2+y$ is not. 
\end{example0}
Recall that a free module is simply a number of copies of the underlying ring.  In the graded setting an important type of free module involves a {\em shift} in grading.
\begin{example0} Define $R(-i)$ as a copy of $R$ in which the unit $1$ appears in degree $i$. Hence $R(-i)_j = R_{j-i}$. Continuing with $R=\RR[x,y]$, we have $R(-2)_1 =0$, $R(-2)_{2} \simeq \RR$, and $R(-2)_3 = R_1 \simeq \RR x \oplus \RR y \simeq \RR^2$.
\end{example0}
\begin{definition} If $M$ is a graded module then the  Hilbert function is $HF(M,d) = \dim_{\RR}M_d$ and the Hilbert series is $HS(M,t) = \sum_{i \in \ZZ}\dim_{\RR}M_it^i.$ 
\end{definition}
For a finitely generated graded $R$-module $M$, the Hilbert function becomes polynomial for $d\gg 0$ \cite[Theorem 2.3.3]{sBook} so is denoted $HP(M,d)$. One raison d'etre for graded objects is that they provide a very clean way to approach the dimension question, as illustrated by the following example.
\begin{example0}\label{BilleraRoseMatrix}
Suppose $\hat \Delta$ is the simplicial complex obtained by embedding $\Delta$ in the plane $\{z_{k+1}=1\} \subseteq \mathbb{R}^{k+1}$
and forming the cone with the origin.  In \cite{br1}, Billera and Rose observe that $S^r(\hat\Delta)$, namely the set of 
splines (of all degrees) on $\hat \Delta$, is a graded module 
over $\RR[x_1,\ldots,x_{k+1}]$ and moreover
$S^r(\hat\Delta)_d \simeq S^r_d(\Delta).$ 

This means computing $\dim_{\RR} \Srd$ is equivalent to computing $HF(\Srhat,d)$, and that for $d \gg 0$ the dimension $\dim_{\RR} \Srd$ is given by the Hilbert polynomial of $S^r(\hat\Delta)$.  When $\Delta \subseteq \RR^k$ is a $k$--dimensional simplicial or polyhedral complex, Billera and Rose note that there is a graded exact sequence:
\[
0\longrightarrow \Srhat \longrightarrow R^{f_k}\oplus
R^{f_{k-1}^0}(-r-1) \stackrel{\phi}{\longrightarrow}
R^{f_{k-1}^0}
\longrightarrow N \longrightarrow
0,
\]
\begin{equation}\label{BRmatrix}
 \hbox{where the map } \phi = \;\;{\small \left[ \partial_k \Biggm| \begin{array}{*{3}c}
l_{\tau_1}^{r+1} & \  & \  \\
\ & \ddots & \  \\
\ & \ & l_{\tau_m}^{r+1}
\end{array} \right]}.
\end{equation}
\end{example0}
Here $f_i$ denotes the number of faces of $\Delta$ of dimension $i$ and $f_i^0$ denotes the number of such faces in the interior of $\Delta$. Notice that this is exactly the matrix appearing in Example~\ref{ex:splineConditionsEx1}. For reasons that will become clear later in this section, we write $[\partial_k \mid D]$ for $\phi$, as well as $N=\coker(\phi)$ and $\Srhat = \ker(\phi)$.

To describe $\partial_k$, note 
that the rows of $\partial_k$ are indexed by $\tau \in \Delta_{k-1}^0$. 
If $\sigma_1, \sigma_2$ denote the $k-$faces adjacent to $\tau$ 
then the smoothness condition means that the only nonzero entries in the row corresponding to
$\tau$ 
occur in the columns corresponding to $\sigma_1$ and $\sigma_2$, and the entries are
$\pm(+1,-1)$. 

\subsubsection{Simplicial complexes}
Homology has its roots in the desire to construct algebraic invariants for topological spaces. To do this, we first create a set of lego blocks (simplices), and then stipulate how the lego blocks can be glued together, yielding a {\em simplicial complex}. Then from the simplicial complex, we construct a sequence of algebraic objects (a {\em chain complex}), which captures topological features. Fasten your seatbelt, here we go!
\begin{definition}\label{simplicialComplex}
An abstract $n$-simplex is a set consisting of all subsets of an 
$n+1$ element ground set. Typically a simplex is viewed as a geometric
object; for example a two-simplex on the set $\{a,b,c\}$ can be visualized
as a triangle, with the subset $\{a,b,c\}$ corresponding to the whole 
triangle, $\{a,b\}$ an edge, and $\{a\}$ a vertex. 

A simplicial complex $\Delta$ on a vertex set $V$ is a collection of subsets
$\rho$ of $V$, such that if $\rho \in \Delta$ and $\tau \subseteq \rho$, 
then $\tau \in \Delta$. If $|\rho| = i+1$ then $\rho$ is called an $i-$face. We write $\Delta_i$ for the set of all $i$-faces of $\Delta$ and $\Delta_i^0$ for the interior $i$-faces. 

An oriented simplex is a simplex with a fixed ordering of the vertices,
modulo the following equivalence relation: for each permutation $\sigma \in S_n$ and
oriented simplex $\tau =[i_1,\ldots, i_n]$, we have $\tau \sim (-1)^{\scriptsize\mbox{sgn}(\sigma)} \sigma(\tau)$.
\end{definition}
\subsubsection{Homology basics and the Euler characteristic}

A sequence of vector spaces and linear transformations
$$\cc:\mbox{ }\cdots\xrightarrow{\; \phi_{j+2}\;} V_{j+1} \xrightarrow{\; \phi_{j+1}\;} V_{j} 
\xrightarrow{\; \phi_{j}\;} V_{j-1} 
\xrightarrow{\; \phi_{j-1}\;} \cdots$$
is called a {\em chain complex} if for all $j$ 
$$\im(\phi_{j+1}) \subseteq \ker(\phi_{j}).$$
The sequence is {\em exact} at position $j$ if 
$\im(\phi_{j+1}) = \ker(\phi_{j})$.  A complex
which is exact everywhere is called an {\em exact sequence}. 
\begin{definition}\label{homologyDef}
The $j$-th homology of the complex $\cc$ is 
\[
H_j(\cc) = \ker(\phi_{j})/ \im(\phi_{j+1}).
\]
The homology of the complex is the graded module $H_\bullet(\cc)=\bigoplus_j H_j(\cc)$.
\end{definition}
We now connect the simplicial complexes constructed in the previous section with the chain complex defined above 
\begin{definition}\label{simpComplexchainComplex}
Let $R$ be a ring and let $C_i(\Delta)$ be the free $R$-module with basis indexed by the oriented $i$-simplices. Elements of $C_i(\Delta)$ are often called the {\em i-chains}. Define a map $C_i(\Delta)\stackrel{\partial_i}{\longrightarrow}C_{i-1}(\Delta)$ 
via
\vskip -.15in
\[
\partial_i[e_{j_0},\ldots ,e_{j_i}] =\sum\limits_{m=0}^i (-1)^m  [e_{j_0},\ldots, \widehat{e_{j_m}}, \ldots ,e_{j_i}]
\]
Here $\widehat{e_{j_m}}$ means to omit $e_{j_m}$. An easy check shows that $\partial_i \circ \partial_{i+1} =0$ and so $\cc = (C(\Delta), \partial)$ is a chain complex.  We define $H_i(\cc)$ as the $i$-th homology of $\cc$. 
\end{definition}
Let's put this into practice on an example.
\begin{example0}\label{orient1}
Consider a simplicial complex with vertices labelled $\{0,1,2\}$ and oriented edges $\{[0,1],[1,2],[2,0]\}$. This is just a hollow triangle, hence topologically equivalent to the circle $S^1$. 
Applying Definition~\ref{simpComplexchainComplex} with the ring $R = \RR$ yields
\[
\cc: 0 \longrightarrow V_1 {\stackrel{\phi}{\longrightarrow}} V_0 \longrightarrow 0,
\]
where $V_1=V_0=\RR^3$ and $\phi$ is:
$$
\left[ \begin{array}{rrr}
-1&0&1\\
1&-1&0\\
0 &1&-1
\end{array} \right]. $$
$H_1(\cc) = \ker(\phi)$ has basis $[1,1,1]^t$ and 
$H_0(\cc) = \coker(\phi) = \RR^3/\im(\phi)$. 
The nonvanishing of $H_1$ captures the fact that $S^1$ has a one-dimensional hole. 

Try repeating this exercise with the boundary of a tetrahedron, which has four triangular faces, six edges, and four vertices. You'll find that $H_1 = 0$ in this case, reflecting the fact that any $S^1$ bounds a set of triangular faces so can be pulled across those faces and contracted down to a point.  However, $H_2$ is now one-dimensional, reflecting the fact that $S^2$ has a two-dimensional ``hole''. 
\end{example0}
For a complex of finite dimensional vector spaces
\[
\cc: 0 \longrightarrow V_n \longrightarrow V_{n-1} \longrightarrow \cdots
\longrightarrow V_1\longrightarrow V_0 \longrightarrow 0
\]
the alternating sum of the dimensions is the 
\textit{Euler characteristic} of $\cc$, written $\chi(\cc)$. When $\cc$ is exact 
$\chi(\cc)=0$. An induction as in \cite[\S 2.3]{sBook} shows that in general
\begin{equation}\label{Euler}
\chi(\cc) = \sum_{i=0}^n (-1)^i \dim V_i =\sum_{i=0}^n (-1)^i \dim H_i(\cc).
\end{equation}
All of the above generalizes directly to sequences
of modules and homomorphisms. 

\subsubsection{The Billera chain complex, and a variant}
Equation~\eqref{compatibility} gives necessary and sufficient conditions for a function $F$ on  a simplicial complex $\Delta \subseteq \RR^k$ to be a spline of smoothness $r$. Equivalently, for a pair of top-dimensional faces $\sigma, \sigma' \in \Delta_k$ such that
\[
\sigma \cap \sigma' = \tau \in \Delta_{k-1}
\]
the pair $(F|_{\sigma}, F|_{\sigma'})$ must be in the kernel of the map
\begin{equation}\label{smoothF}
R^2
\xrightarrow{\; [1,-1]\;}
R/l_{\tau}^{r+1}.
\end{equation}

\begin{example0}
In {\em relative} homology, we pick a subcomplex $\Delta' \subseteq \Delta$ and ``ignore'' it. Algebraically, $C_i(\Delta') \subseteq C_i(\Delta)$, and relative homology is the homology of the chain complex whose $i^{th}$ term is $C_i(\Delta)/C_i(\Delta')$; the differential on $C_i(\Delta)$ induces a differential on the quotient. We have not imposed boundary constraints on the splines studied in this note, so we are computing relative homology, with $\Delta'= \partial(\Delta)$. For Example~\ref{ex:first}, applying
Definition~\ref{simpComplexchainComplex} for the relative complex yields
\[
\partial_2 =\left[ \begin{array}{rrrr}
1&-1&0&0\\
0&1&-1&0\\
0&0 &1&-1\\
-1&0&0 &1
\end{array} \right]
\]
which is the $4 \times 4$ submatrix of $\phi$ from  Example~\ref{ex:splineConditionsEx1}.  We compute 
\[
\partial_2(f_1,f_2,f_3,f_4)^t =\left( \begin{array}{c}
f_1-f_2\\
f_2-f_3\\
f_3-f_4\\
f_4-f_1
\end{array} \right). 
\]
This is the right hand side of Equation~\eqref{firstE}, just as we asserted in Example~\ref{ex:splineConditionsEx1}. 
\end{example0}
We still need to encode the smoothness condition \eqref{compatibility}.  Equation~\eqref{smoothF} provides the clue:
rather than having $\partial_2$ map a free module to another free module, we enrich our chain complex to include the smoothness condition as follows. 
\begin{example0}\label{Alltogether}
Continuing with the previous example, define a map
\[
\bigoplus\limits_{\sigma \in \Delta_2}R \stackrel{\partial_2}{\longrightarrow}\bigoplus\limits_{\tau \in \Delta_1^0}R/l_{\tau}^{r+1} \mbox{ where }\Delta_1^0 \mbox{ denotes interior edges of }\Delta.
\]
When $\Delta_1$ consists of the edges in Example~\ref{ex:first}, the kernel of the map $\partial_2$ consists exactly of polynomial vectors $(f_1,f_2,f_3,f_4)$ satisfying Equation~\eqref{firstE}. 
\end{example0}

This suggests we can define a chain complex whose top homology module is the splines on $\Delta$. Billera defined the following complex in \cite{bTAMS}. 

\begin{definition}\label{Billcomplex}
For each vertex $v \in \Delta_0$, let $I_v$ be the 
ideal of $v \in \RR^{2}$, so $I_v$ is
generated by a pair of linear forms. Fix $r \in \NN$ and 
define the following chain complex
\begin{equation}\label{BilleraSES}
{\mathcal R}/{\mathcal I}: 0 \longrightarrow 
\bigoplus\limits_{\sigma \in \Delta_2}R \stackrel{\partial_2}{\longrightarrow}\bigoplus\limits_{\epsilon \in \Delta_1^0}R/l_{\epsilon}^{r+1} \stackrel{\partial_1}{\longrightarrow}\bigoplus\limits_{v \in \Delta_0^0}R/I_{v}^{r+1} \longrightarrow 0.
\end{equation}
\end{definition}
It is not hard to show that this is a complex, and by construction $\Srhat = H_2({\mathcal R}/{\mathcal I}).$ To prove Strang's conjecture on the dimension of $S^1_d(\Delta)$, Billera used rigidity results of Whiteley \cite{w} to show that when $r=1$ and $\Delta$ is a \textit{generic} triangulation, 
\[
H_1({\mathcal R}/{\mathcal I}) = 0 = H_0({\mathcal R}/{\mathcal I}).
\]
The result follows by computing the Euler characteristic of Equation~\ref{BilleraSES}. What about other $r$ and $\Delta$? It is easy to show that if $\Delta$ is connected then $H_0({\mathcal R}/{\mathcal I})=0$ for any $r$.  However if $r>1$ then $H_1({\mathcal R}/{\mathcal I})$ need not vanish even when $\Delta$ generic. \cite{ss97a} studies the dimension of $\Srd$ by modifying the complex ${\mathcal R}/{\mathcal I}$ in the last position. 
\begin{definition}\label{ssSES}
For a vertex $v \in \Delta$ incident to edges $\epsilon_1, \ldots, \epsilon_n$, define 
\[
J_v  = \langle l_{\epsilon_1}^{r+1}, \ldots, l_{\epsilon_n}^{r+1} \rangle.
\]
The chain complex ${\mathcal{R}}/{\mathcal J}$ replaces $\bigoplus\limits_{v \in \Delta_0^0}R/I_{v}^{r+1}$ with $\bigoplus\limits_{v \in \Delta_0^0}R/J_{v}$ in Equation~\eqref{BilleraSES}.
\end{definition}
This modified complex pays dividends for planar splines \cite{ss97a, ss97b}, but the bigger benefit occurs in higher dimensions \cite{s1}, where in general the complexes ${\mathcal R}/{\mathcal I}$ and ${\mathcal R}/{\mathcal J}$ only agree at the top two terms. 
\subsection{Computational Methods}
All of our computations can be done using the {\tt AlgebraicSplines} package in the free software system {\tt Macaulay2} \cite{danmike}.  The script takes the vertex locations $V$, the maximal simplices $F$, and the order of smoothness $r$. It allows computation of both the module of splines and the chain complex ${\mathcal R}/{\mathcal J}$ of the previous section:
\begin{example0}\label{codeExample1}
We illustrate how to compute the dimension of the spline spaces for the triangulation $\Delta_1$ in Example~\ref{ex:first}, namely the second column of Table 1. 
\begin{verbatim}
i1 : loadPackage "AlgebraicSplines"
 
i2 : V={{0,0},{0,1},{1,0},{0,-1},{-1,0}}

o2 = {{0, 0}, {0, 1}, {1, 0}, {0, -1}, {-1, 0}}

o2 : List

i3 : F={{0,1,2},{0,2,3},{0,3,4},{0,4,1}}

o3 = {{0, 1, 2}, {0, 2, 3}, {0, 3, 4}, {0, 4, 1}}

o3 : List

i4 : splineModule(V,F,0)

o4 = image | 1 t_0 0   t_0t_1 |
           | 1 t_0 t_1 0      |
           | 1 0   t_1 0      |
           | 1 0   0   0      |


--Notice that in this simple example we can read off the splines.
--The first column represents splines for which the same 
--polynomial is assigned to each triangle, the second and third 
--columns reflect the two symmetries: take the linear polynomial 
--t_0 on both right hand sides of the vertical axis, and zero on 
--the other side. Similary, take t_1 above the horizontal axis,
--and zero below it. Finally, take t_0t_1 in the first orthant, 
--and zero elsewhere.

--Next, we find the Hilbert polynomial, where r goes from 0 to 4.

i5 : scan(5,i->(print hilbertPolynomial
                    (splineModule(V,F,i), Projective=>false)))
  2
2i  + 2i + 1
  2
2i  - 2i + 4
  2
2i  - 6i + 13
  2
2i  - 10i + 28
  2
2i  - 14i + 49
\end{verbatim}
\end{example0}
\begin{example0}\label{codeExample2}
Now consider the two triangulations $\Delta_4$ and $\Delta_5$ appearing in Example~\ref{ex:second}. 
\begin{verbatim}
--Adding a semicolon at the end of a line suppresses the output. 
--Because the spline module in the example below has a very 
--complicated presentation, we don't want to see it. We also don't
--need to see the vertex locations and maximal faces echoed back, 
--so we also add a semicolon to suppress those echos.

i2 : T5verts ={{0,8},{2,-2},{0,2},{-2,-2},{8,-6},{-8,-6}};

i3 : T4verts ={{0,5},{1,0},{-1,0},{0,-3},{6,-4},{-6,-4}};

i4 : TriEx2 = {{0,1,2},{0,1,4},{0,2,5},{1,2,3},{2,3,5},{1,3,4},{3,4,5}};

i5 : Spline4r1 = splineModule(T4verts, TriEx2,1); 

--Compute the spline module for Delta_4, r=1

i6 : Spline5r1 = splineModule(T5verts, TriEx2,1); 

--Compute the spline module for Delta_5, r=1

i7 : hilbertPolynomial(Spline4r1, Projective=>false)

     7 2   15
o7 = -i  - --i + 7
     2      2

i8 : hilbertPolynomial(Spline5r1, Projective=>false)

     7 2   15
o8 = -i  - --i + 7
     2      2

--as expected, the hilbertPolynomials agree. 
--What about the Hilbert function?

i9 : hilbertFunction(2,Spline4r1)

o9 = 7

i10 : hilbertFunction(2,Spline5r1)

o10 = 6

--As Morgan and Scott discovered, dim S^1_2 is different for 
--triangulations Delta_4 and Delta_5. What about degree 3?

i11 : hilbertFunction(3,Spline4r1)

o11 = 16

i12 : hilbertFunction(3,Spline5r1)

o12 = 16

--Last, we build the chain complex which computes the splines.

i13 : C5=splineComplex(T5verts,TriEx2,1);

i14 : prune HH_1 C5

o14 = 0

--HH_1 computes H_1 homology; "prune" simplifies the presentation.


i15 : C4=splineComplex(T4verts,TriEx2,1);

i16 : prune HH_1 C4

o16 = cokernel {2} | t_2 t_1 t_0 |

i24 : scan(5, i->print hilbertFunction(i,o16))

0
0
1
0
0

--For the complex Delta_5, the first homology H_1 is zero. For the 
--complex Delta_4, H_1 is a one-dimensional vector space, which
--appears in degree two. This is because the cokernel of the map 
--R^3-->R given by | t_2 t_1 t_0 | (where R is a polynomial ring 
--in 3 variables) is clearly just a one-dimensional vector space. 
--It is in fact a graded vector space; the {2} indicates that the 
--space is generated in degree 2. This is exactly the discrepancy 
--between the Hilbert function and the Hilbert polynomial of the 
--spline space on the Morgan-Scott triangulation.
\end{verbatim}
\end{example0}
\subsection{Splines meet Equivariant Cohomology: the Stanley-Reisner ring}

Let $\Delta$ be a simplicial complex in $\RR^k$. In \cite{c}, Courant introduced a family of piecewise linear continuous functions on $\Delta$ that we now define.
\begin{definition}\label{CourantFunction}
Let $\Delta \subseteq \RR^k$ be a simplicial complex.  The {\em star} of a vertex $v$ is the subcomplex of $\Delta$ consisting of all faces that contain $v$, and is denoted $\stv$. If $v \in \Delta$ is a vertex, the Courant function $\delta_v$ is the unique continuous, piecewise linear function which is zero outside of $\stv$, zero on the boundary of $\stv$, and one at $v$.
\end{definition}
The Courant function $\delta_v$ is easy to visualize: it is a little tent, with centerpole of height one positioned at vertex $v$, and a tent peg staking down the tent at each vertex of the boundary of $\stv$. Recall that  the set $S^0(\Delta)$ of continuous splines has the extra structure of a \textit{ring} because multiplying continuous splines preserves continuity.
\begin{example0}\label{ex:third}
For example, all three simplicial complexes in Example~\ref{ex:first} are combinatorially equivalent, so the Courant functions will behave in the same way. For $\Delta_1$, label the center vertex $v_0$ and label the remaining vertices $v_1,\ldots v_4$, starting with $v_1$ at the top of the figure and moving clockwise. Since $S^0(\Delta)$  is a ring, we may bundle all the vector spaces $S^0_d(\Delta_1)$ into a single structure rather than breaking out splines by degree.  We do this below:

\begin{center}
\vskip -.3in    \[
    \begin{array}{ccccc}
S^0(\Delta_1)&=&
\bigcup\limits_{d=0}^{\infty} S^0_d(\Delta_1) 
 & \simeq &
 \RR[\delta_{v_0},\ldots,\delta_{v_4}]/\Big\langle \sum\limits_{i=0}^4 
 \delta_{v_i} - 1,\delta_{v_1}\delta_{v_3},\delta_{v_2}\delta_{v_4} \Big\rangle.
 \end{array}
\]
\end{center}
To see why each relation holds, first note that the sum of all the Courant functions is a function whose value is one at every point of $\Delta_1$, yielding $\sum \delta_{v_i}=1$. The vanishing of the term $\delta_{v_1}\delta_{v_3}$ reflects the fact that the tents with centerpoles at $v_1$ and $v_3$ do not overlap and so the corresponding functions multiply to zero (respectively  $\delta_{v_2}\delta_{v_4}$ and the vertices $v_2$ and $v_4$). We will see in Theorem~\ref{SRisSpline} that these are all the relations.
\end{example0}

\noindent One theme of this section is that it can be fruitful to view discrete or combinatorial structures through an algebraic lens. In the case of a simplicial complex, one such lens is the {\em Stanley-Reisner ring}.  (We've focused on simplicial complexes embedded in a Euclidean space but an abstract simplicial complex need not be.)
\begin{definition}\label{SRring}
Let $\Delta $ be a simplicial complex on the vertex set $\{v_1,\ldots, v_n\}$ and let $R$ be the polynomial ring $\RR[x_1,\ldots, x_n]$. The Stanley-Reisner ring $SR_\Delta$ is $R/I_\Delta$,
where
\[
  I_\Delta = \langle x_{i_1}\cdots x_{i_j} \mid [v_{i_1},\ldots, v_{i_j}]
  \mbox{ is not a face of }\Delta \rangle.
\]
\end{definition}
The ideal $I_\Delta$ encodes all the non-faces of $\Delta$, so in
particular the simplicial complex $\Delta$ and the ideal $I_\Delta$
carry the same information. In \cite{bADV}, Billera proved the following.
\begin{theorem}\label{SRisSpline}$\cite{bADV}$
For a simplicial complex $\Delta \subseteq \RR^k$ on the vertex set $\{v_1, \ldots v_n\}$, 
\[
SR_\Delta/\Big\langle \sum_{i=1}^n x_i -1 \Big\rangle \simeq \bigcup\limits_{d=0}^{\infty} S^0_d(\Delta) \simeq S^0(\Delta).
\]
\end{theorem}

\noindent To see that there is an inclusion \vskip -.1in
\[SR_\Delta/\Big\langle \sum_{i=1}^n x_i -1 \Big\rangle \hookrightarrow S^0(\Delta),
\]
notice that if $v_{i_0},\ldots, v_{i_m}$ are vertices defining a non-face, then the corresponding Courant function $\delta_{v_{i_0}}\cdots \delta_{v_{i_m}}$ is zero, because the intersection of the supports of the $\delta_{v_{i_j}}$ is empty. This gives the relation from the previous example:
\[
\sum_{v \in \Delta_0}\delta_v = 1
\]

But where is the connection to equivariant cohomology? One answer is that an embedded simplicial complex $\Delta$ can be used to define an object $X_\Delta$ known as a {\em toric variety}. Toric varieties are key examples in algebraic geometry that have also played an essential role in solving several famous problems on convex polytopes; for more on this story, see \cite{CLS}. The bridge is provided by the following theorem.
\begin{theorem}$[\mbox{Bifet-DeConcini-Procesi} \cite{BDP}, \mbox{Brion} \cite{brion}]$
    The equivariant cohomology ring of the simplicial toric variety $X_\Delta$ is the Stanley-Reisner ring of $\Delta$.
    \end{theorem}
This is absolutely remarkable! A fundamental object in approximation theory---the set of piecewise linear continuous functions on a triangulation $\Delta$---is (up to quotient by the innocuous linear form $\sum x_i -1$) the equivariant cohomology ring of a fundamental object in algebraic geometry---the toric variety associated to $\Delta$. When $\Delta$ is polyhedral rather than simplicial, Payne shows the relation to splines remains \cite{payne}. This is a perfect segue to \S \ref{sec:groupActions}, where we flesh out the details of group actions. 
\section{Group actions}\label{sec:groupActions}
In many areas of mathematics, a general theory can be enhanced if there are symmetries. Symmetries often signal the presence of a group action, and in this section, we pursue this in the setting of splines.
 \begin{definition}
 Let $G=(V,E, \ell:E\rightarrow \mathbb{R}[x_1, \ldots, x_n])$ be a graph with edge-labeling. An automorphism $\varphi$ of $G$ is a pair $(\varphi_V,\varphi_\ell)$ where
\begin{itemize}
\item[] \begin{itemize}
    \item[$\bullet$] $\varphi_V:V\rightarrow V$ is a bijection such that $v_1v_2\in E$ if and only if $\varphi_V(v_1)\varphi_V(v_2)\in E$,
   \item[$\bullet$] $\varphi_\ell:\mathbb{R}[x_1, \ldots, x_n]\rightarrow\mathbb{R}[x_1, \ldots, x_n]$ is an invertible map of $\mathbb{R}$-algebras that acts on each edge-label by $\ell(\varphi_V(v_1)\varphi_V(v_2))=\lambda_{v_1,v_2}\varphi_\ell(\ell(v_1v_2))$ for some non-zero $\lambda_{v_1,v_2}\in\mathbb{R}$.
 \end{itemize} 
 \end{itemize}
 \end{definition}
 Notice that in our case, an invertible map of algebras just means that $\varphi_\ell$ performs a change of basis on the linear span of $\{x_1, \ldots, x_n\}$.

We can of course compose (and invert) maps from a set to itself, so that the set of automorphisms $Aut(G)$ of an edge-labelled graph $G$ is equipped with a group structure, where the operation is just given by map composition on each component: 
\[\varphi\circ \psi=(\varphi_V\circ \psi_V, \varphi_\ell\circ\psi_\ell) \quad \textup{ for any pair } \quad \varphi=(\varphi_V,\varphi_\ell), \psi=(\psi_V, \psi_\ell).\]

 \begin{example0}\label{ex:AutomorphismDelta1}
 The triangulation $\Delta_1$ in Figure 1 (see Example~\ref{ex:first}) and its corresponding  dual edge-labeled graph $G$ in Figure \ref{figure : graph of Delta 1} are very symmetric. The symmetry yields an automorphism $\varphi$ of $G$ that is defined by 
 \[
 \varphi_V\ :\ v_1\mapsto v_3, \ v_2\mapsto v_2, \ v_3\mapsto v_1, \ v_4\mapsto v_4, \qquad 
 \varphi_\ell \ : \ x\mapsto y, \ y\mapsto x.\]
Composing $\varphi$ with itself gives the trivial automorphism, so $\varphi$ has order two. 
 \end{example0}
 
It is intuitively clear that there are rigid motions of the plane that preserve the triangulation $\Delta_1$ of Example~\ref{ex:AutomorphismDelta1}.  Whenever that is the case, the dual edge-labeled graph will inherit corresponding automorphisms.  But there are automorphisms of edge-labeled graphs that do not come from rigid motions of the plane, providing a larger family of maps.  Even better, these automorphisms induce a map on splines that gives us an additional tool with which to analyze splines.
 
 Let $G$ be an edge-labeled graph and let $\varphi=(\varphi_V, \varphi_\ell)$ be an automorphism of $G$. Observe that $\varphi$ induces an automorphism of $S^r_d(G)$ considered as an $\mathbb{R}$-algebra:
 \begin{equation}\label{eqn:ActionOnSpline}
\varphi\cdot f_\bullet=f'_\bullet   
 \quad \hbox{where }\quad f'_v=\varphi_{\ell}(f_{\varphi_V^{-1}(v)}).
 \end{equation}
To see this is well defined, say $f_\bullet\in S^r(G)$ and let $(v_1,v_2)$ be an edge whose preimage is the edge $(w_1,w_2):=(\varphi_V^{-1}(v_1), \varphi_V^{-1}(v_2))$.  Since $f_\bullet$ is a spline we have $\ell(w_1 w_2)^{r+1}\mid f_{w_1}-f_{w_2}$.  Thus if we write 
\begin{align*}
 f'_{v_1}-f'_{v_2}&=\varphi_\ell\left(f_{\varphi_V^{-1}(v_1)}\right)-\varphi_\ell\left(f_{\varphi_V^{-1}(v_2)}\right)=\varphi_\ell\left(f_{w_1}-f_{w_2}\right),   
\end{align*}
 the last term is divisible by $\varphi_\ell(\ell(w_1w_2)^{r+1})$.  By the definition of automorphisms, this is a (non-zero) multiple of $\ell(v_1 v_2)^{r+1}$. 
 
 This shows that an automorphism of an edge-labelled graph induces an automorphism of the corresponding spline space. It is immediate to see that the identity automorphism of $G$ induces the identity automorphism on $S^r(G)$ and that this process respects composition, in the sense that $\varphi\circ \psi$ induces the same automorphism on $S^r(G)$ as the composition of the automorphisms induced by $\varphi$ and $\psi$. In short, the space of splines $S^r(G)$ is equipped with the structure of an $Aut(G)$-representation.
 
More generally, if $H \leq Aut(G)$ is a subgroup of the automorphism group of $G$, then letting every element of $H$ act via Equation~\eqref{eqn:ActionOnSpline} gives splines $S^r(G)$ the structure of an $H$-representation. We observe that every automorphism preserves the degrees of the elements in $S_r(G)$, so for every integer $d \geq 0$ the subspace $S^r_d(G)$ also has the structure of an $H$-representation.

Whenever we are given a group representation, that is, a vector space $U$ equipped with the group action, it is natural to try to decompose the representation into its smallest building blocks.  

\begin{definition}
    Let $U$ be vector space equipped with a group action. A proper non-zero subspace $W$ of $U$ is a (proper) \textit{subrepresentation} if it is  preserved by the group action  (globally, not pointwise).  A representation $U$ is \textit{simple} if it contains no proper subrepresentation.   $U$ is a \textit{trivial representation} if $g \cdot u = u$ for all $g \in G$ and $u \in U$.
    \end{definition}

We noted earlier that $S^r_d(G)$ is a subrepresentation of $S^r(G)$.  The constant splines $M = \{ p(1,1, \ldots, 1)=(p,p,\ldots,p): p \in \mathbb{R}[x_1,\ldots,x_n]\}$ are another example of a subrepresentation. Acting by $\varphi \in Aut(G)$ on $f_\bullet = p(1,1,\ldots,1)$ gives $\varphi(f_\bullet) \in M$ since
\begin{equation} \label{equation: constant splines}\varphi \cdot f_\bullet =  (\varphi_\ell(p), \varphi_\ell(p),\ldots,\varphi_\ell(p)) = \varphi_\ell(p)(1,1,\ldots,1). \end{equation}

Simple representations are harder to identify, but one-dimensional representations are always simple since their only proper (vector) subspace is zero.

All splines in a trivial representation are common eigenvectors to all group elements with eigenvalue $1$.  Constant splines play a special role in trivial representations.  For example, suppose $f_\bullet \in S_0^r(G)$ is both constant and degree zero, namely $p=\mu\in\mathbb{R}$ in the definition of $M$ above.   All $\mathbb{R}$-algebra homomorphisms are by definition $\mathbb{R}$-linear so $\varphi_\ell(\mu)=\mu$ in Equation~\eqref{equation: constant splines} and hence $\varphi \cdot f_\bullet = f_\bullet$ for all $\varphi \in Aut(G)$.

 \begin{example0}
 Consider the subgroup $H$ of $Aut(G)$ generated by $\varphi$ from Example \ref{ex:AutomorphismDelta1}. We want to consider the $H$-action on $S^r(G)$, which is of course uniquely determined by the action of $\varphi$ on $S^r(G)$. 

Since $\varphi$ has order 2, $H$ is isomorphic to a cyclic group of order two. However $\varphi$ does not act as multiplication by $\pm 1$ on the spline space.  Indeed,
the element $f_\bullet=(x,x,x,x)$ belongs to $S^r(G)$ for any $r$ but $\varphi\cdot f_\bullet=(y,y,y,y)$, which is of course not equal to $\pm(x,x,x,x)$. Hence we deduce that the real vector space $S^r(G)$ decomposes into a direct sum of two non-zero eigenspaces $V(\varphi,1)$ and $V(\varphi,-1)$ relative to the eigenvalues $1$ and $-1$. Of course, since $\varphi$ preserves degree, we can restrict to any degree $d$ to get $S^r_d(G)=V(\varphi,1)_d\oplus V(\varphi,-1)_d$, where $V(\varphi, i)_d$ denotes the eigenspace of eigenvalue $i$ in $S^r_d(G)$. However, it might happen that one $V(\varphi,i)_d$ is zero --- as in the case of $d=0$.

In fact, determining the dimensions of the eigenspaces $V(\varphi,1)_d$  and $V(\varphi,-1)_d$ is equivalent to decomposing $S_d^r(G)$ into irreducible $H$-representations. More generally, the cyclic group of order two has only two isomorphism classes of irreducible representations: the trivial representation $\textrm{Triv}$ (1-dimensional representation on which the group acts trivially) and the sign representation $\textrm{Sgn}$ (1-dimensional representation on which the non-neutral element acts by multiplication by $-1$). Therefore, as an $H$- representation we have
\[
S_d^r(G)\simeq \underbrace{\textrm{Triv}\oplus  \ldots\oplus \textrm{Triv}}_{m_d-times}\oplus \underbrace{\textrm{Sgn}\oplus \ldots\oplus\textrm{Sgn}}_{n_d-times},
\]
where $m_d=\dim V(\varphi,1)_d$ and $n_d=\dim V(\varphi,-1)_d$.


For example, if $r=0$ and $d=1$, then $S^0_1(G)$ has an $\mathbb{R}$-basis given by
\[
(1,1,1,1), \ (x,x,x,x), \ (y,y,y,y), \ (x,x,0,0), \ (0,y,y,0).
\]
We can adjust it to a basis that is compatible with the eigenspace decomposition:
\[
\underbrace{(1,1,1,1), (x+y,x+y,x+y,x+y), (x,x+y,y,0)}_{\text{have eigenvalue 1}},\]
\[\underbrace{(x-y,x-y,x-y,x-y), \ (x,x-y,-y,0)}_{\text{have eigenvalue -1}}.
\]
Therefore $S^0_1(G)\simeq \textrm{Triv}\oplus \textrm{Triv}\oplus \textrm{Triv} \oplus\textrm{Sgn}\oplus \textrm{Sgn}$.
 \end{example0}

Observe that symmetries also give us information about the $\mathbb{R}[x,y]$-module structure of the spline module.  In the above example, before looking at a basis for $S^0_1(G)$ we could have noticed that the constant splines $M \cap S^0_1(G)$ form a subrepresentation of $S^0_1(G)$.  In addition, the nonconstant spline $f_\bullet=(x,x,0,0)$ is not sent by $\varphi$ to a scalar multiple of itself.  Hence $S^0_1(G)$  contains the $\mathbb{R}$-span of $f_\bullet$ and $\varphi\cdot f_\bullet$ together with $M \cap S^0_1(G)$, so has at least dimension five.  
\pagebreak

\begin{example0}\label{TriTriExample}
Consider the polyhedral subdivision $\Delta_4$ depicted in Figure~\ref{figure : TriTri} below.
\begin{figure}[h]
\begin{center}
\epsfig{file=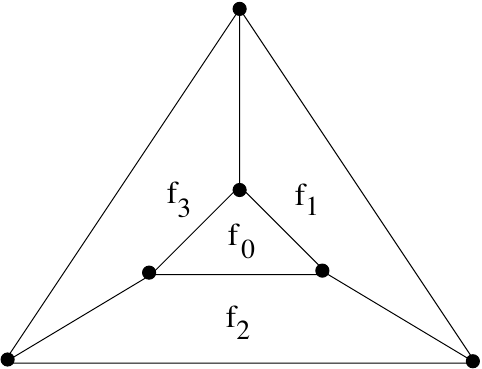,height=1.3in,width=1.5in}
\end{center}
\caption{Polyhedral subdivision of a triangle.}
\label{figure : TriTri}
\end{figure}
If both triangles are equilateral and centered at the origin, and we assume that the top vertex of the interior triangle is $(0,1)$, then the corresponding edge-labeled graph is in Figure \ref{figure : graph polyhedral subdivision}.
\begin{figure}
    \centering
\begin{tikzpicture}[x=0.75pt,y=0.75pt,yscale=-1,xscale=1]

\draw   (247.3,23.56) -- (415.27,23.56) -- (331.29,169.02) -- cycle ;
\draw    (247.3,23.56) -- (331.29,72.05) ;
\draw    (331.29,72.05) -- (331.29,169.02) ;
\draw    (331.29,72.05) -- (415.27,23.56) ;
\draw (368.49,112.03) node [anchor=north west][inner sep=0.75pt]  [font=\scriptsize,rotate=-301.25] [align=left] {$\displaystyle x+\sqrt{3} y$};
\draw  (333,75.20) node [anchor=north west][inner sep=0.75pt]  [font=\scriptsize,rotate=0] [align=left] {$v_0$};
\draw  [fill={rgb, 255:red, 0; green, 0; blue, 0 }  ,fill opacity=1 ](331.29,75.20) .. controls (329.8,74.86) and (328.6,73.6) .. (328.6,72.05) .. controls (328.6,70.49) and (329.8,69.24) .. (331.29,69.24) .. controls (332.77,69.24) and (333.98,70.49) .. (333.98,72.05) .. controls (333.98,73.6) and (332.77,74.86) .. (331.29,74.86) -- cycle ;
\draw  (232.3,20.37) node [anchor=north west][inner sep=0.75pt]  [font=\scriptsize,rotate=0] [align=left] {$v_3$};
\draw  [fill={rgb, 255:red, 0; green, 0; blue, 0 }  ,fill opacity=1 ] (247.3,26.37) .. controls (245.82,26.37) and (244.61,25.11) .. (244.61,23.56) .. controls (244.61,22.01) and (245.82,20.75) .. (247.3,20.75) .. controls (248.79,20.75) and (249.99,22.01) .. (249.99,23.56) .. controls (249.99,25.11) and (248.79,26.37) .. (247.3,26.37) -- cycle ;
\draw (421.3,20.37) node [anchor=north west][inner sep=0.75pt]  [font=\scriptsize,rotate=0] [align=left] {$v_1$};
\draw  [fill={rgb, 255:red, 0; green, 0; blue, 0 }  ,fill opacity=1 ] (415.27,26.37) .. controls (413.79,26.37) and (412.58,25.11) .. (412.58,23.56) .. controls (412.58,22.01) and (413.79,20.75) .. (415.27,20.75) .. controls (416.76,20.75) and (417.96,22.01) .. (417.96,23.56) .. controls (417.96,25.11) and (416.76,26.37) .. (415.27,26.37) -- cycle ;
\draw (329.5,173.84) node [anchor=north west][inner sep=0.75pt]  [font=\scriptsize,rotate=0] [align=left] {$v_2$};
\draw  [fill={rgb, 255:red, 0; green, 0; blue, 0 }  ,fill opacity=1 ] (331.29,171.84) .. controls (329.8,171.84) and (328.6,170.58) .. (328.6,169.02) .. controls (328.6,167.47) and (329.8,166.21) .. (331.29,166.21) .. controls (332.77,166.21) and (333.98,167.47) .. (333.98,169.02) .. controls (333.98,170.58) and (332.77,171.84) .. (331.29,171.84) -- cycle ;

\draw (324,10) node [anchor=north west][inner sep=0.75pt]  [font=\small] [align=left] {$\displaystyle x$};
\draw (286.06,24.57) node [anchor=north west][inner sep=0.75pt]  [font=\scriptsize,rotate=-30.17] [align=left] {$\displaystyle y-\sqrt{3} x-1$};
\draw (335.48,47.43) node [anchor=north west][inner sep=0.75pt]  [font=\scriptsize,rotate=-330.48] [align=left] {$\displaystyle y+\sqrt{3} x-1$};
\draw (333,93) node [anchor=north west][inner sep=0.6pt]  [font=\small] [align=left] {$\displaystyle y+\frac{1}{2}$};
\draw (368.49,112.03) node [anchor=north west][inner sep=0.75pt]  [font=\scriptsize,rotate=-301.25] [align=left] {$\displaystyle x+\sqrt{3} y$};
\draw (268.64,71.2) node [anchor=north west][inner sep=0.75pt]  [font=\scriptsize,rotate=-59.59] [align=left] {$\displaystyle x-\sqrt{3} y$};
\end{tikzpicture}
    \caption{Graph of the polyhedral subdivision in Figure~\ref{figure : TriTri}.}
    \label{figure : graph polyhedral subdivision}
\end{figure}
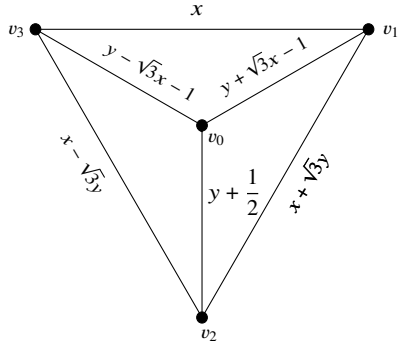
The rotation corresponding to the matrix 
\[
A=\left(
\begin{array}{cc}
    -1/2 &  -\sqrt{3}/2\\
    \sqrt{3}/2 & -1/2 
\end{array}
\right)
\]
sends our triangle and its polyhedral subdivision to itself, hence induces an automorphism $\varphi$ of the corresponding graph, namely:
\[
 \varphi_V\ :v_0\mapsto v_0, \qquad\ v_1\mapsto v_3, \ v_2\mapsto v_1, \ v_3\mapsto v_2, \  \qquad 
 \]
 and 
 \[
 \varphi_\ell \ : \ x\mapsto -\frac{1}{2}x+\frac{\sqrt{3}}{2}y, \quad y\mapsto -\frac{\sqrt{3}}{2}x-\frac{1}{2}y.
 \]
 The automorphism $\varphi$ has order 3 and so it generates a copy of $C_3$, the cyclic group of order 3. Thus we get a $C_3$-representation on $S^r_d(G)$ for any $r$ and $d$.

Consider once again the case when $r=0$  and $d=1$.
We will denote spline elements as a quadruple whose entries agree with the vertex numbering. Let
$$f_\bullet=(f_{v_0},f_{v_1}, f_{v_2}, f_{v_3})=(\sqrt{3}x+y-1, 0,\sqrt{3}x+3y ,2\sqrt{3}x),$$ 
so that 
\[
\varphi\cdot f_\bullet=(y-\sqrt{3}x-1, -2\sqrt{3}x,3y-\sqrt{3}x,0), 
\]
and
\[
\varphi^2\cdot f_\bullet=
(-2y-1,-\sqrt{3}x-3y,0,\sqrt{3}x-3y)
.
\]
 The $\mathbb{R}$-span $U$ of $\{f, \varphi\cdot f, \varphi^2\cdot f\}$ is (globally) sent to itself by $\varphi$ and hence by $C_3$. Thus $U$ is a 3-dimensional sub-$C_3$-representation of $S^0_1(G)$. Moreover, the representation $U$ has itself a subrepresentation of dimension 1 generated by the $\varphi$-invariant element $f+\varphi\cdot f+\varphi^2 \cdot f$.  In other words, we have a copy of the trivial representation in $U$. 
 
 Now look at $\varphi$ as an automorphism of $U$.  If we write $\varphi$ in the basis $\{f,\varphi\cdot f, \varphi^2 f\}$ we get the matrix   $A=\left(
\begin{array}{ccc}
    0 &  0&1\\
1&0&0   \\
0&1&0
\end{array}
\right)$. The eigenvalues of this matrix are the third-roots of unity.  Two of these are complex (non real) eigenvalues, so to find a basis for $U$ of eigenvectors, we have to change the base field.  If we take our coefficients in $\mathbb{C}$, then $U$ splits into a direct sum of three irreducible one-dimensional representations, corresponding to the three eigenvalues. Explicitly, the following three elements of the complexified spline generate three $C_3$-representations:
\[
f+\varphi\cdot f+ \varphi^2\cdot f, \quad f+\xi\left(\varphi\cdot f\right)+ \xi^2\left(\varphi^2\cdot  f\right),
\quad f+\xi^2\left(\varphi\cdot f\right)+ \xi\left(\varphi^2 \cdot f\right)
\]
where $\xi$ is a fixed primitive third root of 1.
\end{example0}

Group representations on splines of geometric origin, such as the ones giving equivariant cohomology of GKM spaces, are particularly interesting to geometric representation theorists, whose more general aim is to construct concrete (geometric) realization of algebraic objects such as representations of groups or algebras.  Finding such a realization means that geometric tools can be used to answer algebraic questions (and vice versa). One advantage of using GKM theory to construct the geometric realization is that one can also exploit the extremely explicit description of the equivariant cohomology, and often the combinatorics of the underlying graph. This has been done, for example, in studying Weyl group representations on the equivariant cohomology of Schubert varieties in \cite{Ty08b}, symmetric group actions on equivariant cohomology of Hessenberg varieties \cite{Ty08}, or on certain varieties of quiver representations \cite{LP21}. All these examples involve only spline spaces for $r=0$. It would be very interesting to study the representations arising on splines for $r\geq 1$, which have not been investigated (so far) by representation theorists.
\section{Open Problems and Further Directions}
There are a myriad of interesting open problems involving splines and connections to combinatorics, algebra, geometry, and topology, and we close this survey with an overview of several such problems; see also the paper of Sottile \cite{Sottile} in this volume, as well as problems in the survey paper \cite{s4}.
\subsection{Open problems connected to algebraic geometry}
\vskip -.1in
We start with an overview of problems with close connections to algebraic geometry. 
\vskip -.1in
\subsubsection{Bounds for planar splines}
It is an open question if a triangulation $\Delta \subseteq \RR^2$ of a simply connected region achieves the lower bound of Equation~\ref{bound} when $r=1$ and $k=3$: is $L^1_3(\Delta) = \dim_{\RR}S^1_3(\Delta)$? Examples of Diener in \cite{diener} and Tohaneanu in \cite{tohaneanu} show that there exist $\Delta$ with $L^r_{2r}(\Delta) \ne \dim_{\RR}S^r_{2r}(\Delta)$, and Alfeld-Schumaker show in \cite{as2} that if $\Delta$ is a generic triangulation, then $L^r_{3r+1}(\Delta) = \dim_{\RR}S^r_{3r+1}(\Delta)$.

What happens in the range $2r+1 \le d < 3r+1$? A conjecture in \cite{sStiller} proposed that $L^r_{2r+1}(\Delta) = \dim_{\RR}S^r_{2r+1}(\Delta)$ for arbitrary planar triangulations. However, in 2019, Yuan \cite{sy} found a triangulation $\Delta_Y$ such that $L^2_{5}(\Delta_Y) \ne \dim_{\RR}S^2_{5}(\Delta_Y)$, and in fact for any $r$, that special triangulation $\Delta_Y$ was shown in \cite{ssy} to satisfy \vskip -.1in
\[
L^r_{\lfloor 2.2r \rfloor}(\Delta_Y) \ne \dim_{\RR}S^r_{\lfloor 2.2r \rfloor}(\Delta_Y).
\]
Hence, the two central open questions on planar splines are 

 \hskip .2in $\bullet$ 
Is $L^1_3(\Delta) = \dim_{\RR}S^1_3(\Delta)$ for all triangulations? 

\hskip .2in $\bullet$ What happens to $\dim_{\RR}\Srd$ in the range $\lfloor 2.2r \rfloor < d <3r+2$?
  \subsubsection{Special configurations: root systems and crosscut partitions} 
In the context of a particular problem, we often want a ``designer" subdivision: a simplicial or polyhedral complex that takes advantage either of the geometry of the region being subdivided, or of the type of PDE/approximation theory problem being solved. In the planar case, the ``crosscut partitions'' of Chui and Wang \cite{cw} are an example, which are generalized to ``pseudoboundary partitions'' in \cite{ss97b}. 

Another example is the {\em Alfeld split} of a simplex. Foucart and Sorokina conjectured a dimension formula for $\Srhat$ in \cite{FoucartT}. For $r=1$ their conjecture was proved in \cite{KolesnikovT}. 

Using results of 
Terao \cite{terao} on hyperplane arrangements and type A root systems, the general conjecture was proved in \cite{s3}. Except for the Alfeld split, there are not many general formulas known for special configurations in dimension three or more; this is a fruitful area for exploration.\subsubsection{Splines in dimension three or more}
As one might guess, given that there are still open questions in dimension two, the situation in higher dimensions becomes even more complicated. Much work has been devoted to {\em tetrahedral splines}, namely for a simplicial subdivision $\Delta$ of a domain in $\RR^3$; see for examples papers by various subsets of Alfeld, DiPasquale, Mourrain, Schumaker, Villamizar, and Whiteley \cite{as3, asw, mdpVillamizar, mv}. 

A cautionary point: in the dimension three case, the vertex ideal $J(v)$ defined in \S 2.2 is an ideal generated by powers of linear forms in three variables. Via the inverse systems approach sketched in \S 4.1.8, this is related to the Segre-Harbourne-Gimigliano-Hirschowitz conjecture on fatpoints on $\pp^2$, which has been open for many years. See Miranda \cite{miranda} for additional details. 

For general $\Delta \subseteq \RR^k$, a rather technical spectral sequence argument in \cite{s1} gives necessary and sufficient conditions for $\Srhat$ to be a free module, giving the conditions in terms of vanishing of lower homology; a more complicated variant of the spectral sequence yields similar results in the polyhedral setting, as shown in \cite{s2}.

\subsubsection{Splines on polyhedral subdivisions}
In applications, we often want to allow non-simplicial subdivisions, for example, polyhedral subdivisions like the one below:
\begin{example0}\label{TH}
Let $P$ be the polygonal complex depicted in Figure~\ref{figure : TriTri} of Example~\ref{TriTriExample} and let $P'$ be a complex obtained by perturbing (any) vertex so that the affine spans of the 
three edges which connect boundary vertices to interior vertices are not concurrent. (In a sense, this is a polyhedral analog of the Morgan-Scott triangulation of Example~\ref{ex:second}.)  The dimensions for $S^r_d(P)$ and $S^r_d(P')$ for small values of $r$ and $d \gg 0$ are below. 
\vskip .1in
\begin{center}
\begin{supertabular}{|c|c|c|}
\hline $r$ & $\dim_{\mathbb{R}} S^r_d(P)$ & $\dim_{\mathbb{R}} S^r_d(P')$ \\
\hline $0$ & $2d^2+2$ & $2d^2+1$ \\ 
\hline $1$ & $2d^2-6d+10$ & $2d^2-6d+7$\\ 
\hline $2$ & $2d^2-12d+32$ & $2d^2-12d+25$\\ 
\hline $3$ & $2d^2-18d+64$ & $2d^2-18d+52$\\ 
\hline $4$ & $2d^2-24d+110$ & $2d^2-24d+91$\\ 
\hline
\end{supertabular}
\end{center}
\end{example0}
\vskip .1in
The values for $S^r_d(P)$ and $S^r_d(P')$ are explained by work of Mcdonald and Schenck in \cite{mcdonald}, which proves a variation of Equation~\ref{bound} for planar polyhedral splines. The key ingredient is to add extra terms $\gamma_i$ that reflect ``ghost'' vertices: in the complex $P$, the three edges connecting boundary vertices to interior vertices are concurrent, whereas in the perturbed complex $P'$ they are not. 

There are many open questions on planar polyhedral splines. For instance, DiPasquale shows in \cite{mdp2} that the Alfeld-Schumaker bound for when the dimension formula becomes a polynomial depends on the maximal number of edges of a polygon in the subdivision.

\subsubsection{Supersmoothness} For any complex $\Delta$, the space of splines $\Srd$ consists of piecewise polynomials which have smoothness {\em at least} $r$ everywhere. But for some special configurations, the condition of smoothness $r$ everywhere in fact results in a {\em higher} order of smoothness across certain faces. This is known as {\em supersmoothness}, first studied by Sorokina in \cite{Sor10, Sor14}. Further work of Shehktman and Sorokina in \cite{borisTanya15, borisTanya}, Floater and Hu in \cite{FloaterHu}, Toshniwal and Villamizar in \cite{TV} has led to a deeper understanding of supersmoothness, but much remains mysterious.

\subsubsection{Splines and standard topological operations} What happens when standard topological operations are applied to $\Delta$? Examples of these operations include 
\begin{itemize}
    \item Subdividing a face (studied in \cite{sTanya}).
    \item Deleting a subcomplex.
    \item The Mayer-Vietoris sequence (\S 4.4.1, \cite{sBook2}), which is an exact sequence obtained from identifying two simplicial complexes $\Delta_1$ and $\Delta_2$ along a common subcomplex $\Delta_1 \cap \Delta_2$.
    \end{itemize}
An application of Mayer-Vietoris appears in \cite{mcdonald2}, but there are a wide range of similar open questions along these lines. \subsubsection{Semialgebraic splines} While the standard subdivisions used are often simplicial or polyhedral, there is increasing interest in partitions where the maximal cells meet along a boundary that is nonlinear. A natural first case to consider is when the boundary is defined by a quadratic or higher-degree polynomial. This situation was first considered by Stiller in \cite{stiller} in the planar case; Stiller studied the question locally and formulated the problem as a question on syzygies, as in \S 2.2. 

Recent work of DiPasquale, Sottile and Sun \cite{mdpSottileSun}, and DiPasquale and Sottile \cite{mdpSottile} has yielded dimension formulas for some very specific situations, and in \cite{borisTanya}, Shekhtman and Sorokina explore a connection to supersmoothness.

\subsubsection{Inverse systems and Fatpoint Ideals}
In Definition~\ref{ssSES} of \S 2.2 we saw that splines appear as the top homology module of a chain complex with terms of the form $R/J$, where $R$ is a polynomial ring and $J$ is an ideal generated by powers of linear forms. 

{\em Apolarity} is an algebraic technique introduced by Macaulay \cite{m}, and work in \cite{ei, i} applies apolarity to translate questions about ideals of powers of linear forms to questions about ideals of {\em fatpoints} in {\em projective space}. This connection is described in detail in \S 5 of \cite{s4}, and is used by Geramita and Schenck in \cite{gs} to give a dimension formula for planar splines of mixed orders of smoothness. 

This line of research ties splines to a famous open question of Segre, Gimigliano, Harbourne and Hirschowitz on fatpoints in the projective plane. See \cite{miranda} for a discussion of their conjecture and the work of Alexander and Hirschowitz \cite{ah1, ah2}.

\subsection{Open problems connected to topology or representation theory}
We now switch our attention to problems with connections to group actions.
\subsubsection{When is the spline module $S^r(G)$ free over the coefficient ring?} 
Free modules ---which we defined in \S 2.2.1--- are the simplest sort of module because they are {\em free} of relations and so behave, in essence, like vector spaces. Contrast this to the module generated by $\{x, y\} \subseteq \RR[x,y]$. There are of course two generators, but also a relation 
\[
y \cdot x -x \cdot y =0.
\]
It turns out that there are cases where the spline module is free. 

The case of most interest for approximation theory occurs when the region $D \subseteq \RR^k$ is topologically trivial, namely $D$ has no ``holes'' of any dimension. If $\Delta$ is a simplicial subdivision of $D$, a spectral sequence argument in \cite{s1} shows that the spline module $S^r(\widehat \Delta)$ is {\em free} if and only if 
\[
H_i(R/J) = 0 \mbox{ for all } i < k,
\]
where $\mathcal{R}/\mathcal{J}$ is again the chain complex from \S 2.2. This means that when the spline module is free, the Hilbert series of $S^r(\widehat \Delta)$ is simply the Euler characteristic of the complex $\mathcal{R}/\mathcal{J}$. In the polyhedral case, the spectral sequence is more complicated, but a similar result holds \cite{s2}. 

In equivariant cohomology, one of the foundational hypotheses of GKM theory is that the torus $T$ acts on $X$ according to a condition called \textit{equivariant formality}. Depending on the mathematical subdiscipline, experts can mean different things by equivariant formality, but one thing all definitions imply is that the equivariant cohomology $H^*_T(X)$ is a free module over the equivariant cohomology of a point $H^*_T(pt)$.  Translating this to our language, alla equivariant cohomology splines $S^r(G_X)$ are assumed to be free over the polynomial ring.

However, when $G$ is an arbitrary edge-labeled graph, then even when the coefficient ring $R$ is a polynomial ring, very little is known in general about when $S^r(G)$ is free over $R$.  When $G$ is a tree then the answer is  yes, as shown in \cite{GTV16}.

\subsubsection{Explicit construction of spline bases}

When splines arise as the equivariant cohomology of some variety, we often have geometric reasons to prefer a particular basis, for instance because a family of subvarieties like Schubert varieties induce a geometric basis, or because the torus action induces a basis through Morse theory or the algebraic-geometric analogue due to Bialynicki-Birula \cite{BB1, BB2}.  In these cases, geometry often provides global information about the basis --- e.g. the degree of each polynomial, or which polynomials are zero --- but very little specific information about each polynomial.

We want to know the specific entries of the basis elements.  For instance, when $X$ is the flag variety, we can represent the elements of the equivariant Schubert basis as splines $p_\bullet \in S^r(G_X)$.  Billey \cite{billey} and Andersen, Jantzen and Soergel \cite{AJS94} separately gave an explicit combinatorial formula for the polynomial $p_v$ at each vertex of each $p_\bullet$.

\subsubsection{Ring structure of splines}

Given an \textit{additive} basis $\{f_{\bullet, i}\}$ for the splines $S^r(G)$, we can expand any product $f_{\bullet, 1} f_{\bullet, 2}$ in terms of the basis as
\[f_{\bullet, 1} f_{\bullet, 2} = \sum_j c_{ij}^k f_{\bullet, k}\]
for some polynomial coefficients $c_{ij}^k$ via the module structure of splines.   
These $c_{ij}^k$ are called \textit{structure constants}.  An important problem in modern Schubert calculus is to find illuminating combinatorial formulas for these structure constants.  

Hermann Schubert was not thinking in terms of ring structure in 1879 when he did the virtuosic linear algebra calculations that gave the field its name; see  \cite{KleimanLaksov} for a survey of how Schubert's data came to be understood as structure constants in the cohomology of the Grassmannian. More recently, Knutson and Tao \cite{KT} gave a formula for these $c_{ij}^k$ in the equivariant cohomology of Grassmannians using the technology of splines, via GKM theory.  Researchers have generalized to Grassmannians of other Lie types, to other cohomology theories, and to other varieties. Many open questions remain, even something as basic as finding structure constants in the ordinary cohomology of the flag variety. For a recent survey of open problems and connections between equivariant and ordinary cohomology, see \cite{RYY}.

\subsubsection{Questions about quotients and cohomology versus equivariant cohomology}

Quotients appear in many parts of our story about splines.  From the dual perspective, it's natural to describe the spline ring $S^r(\Delta)$ as a quotient of a polynomial ring, as in Example~\ref{ex:third}.  When we do this,  what are the relations defining the ring?  In particular, if we describe $S^r(\Delta)$ as a quotient, what goes in the ideal? This is a difficult question, especially if the goal is to find a computationally-useful description of the quotient. 

Another way that quotients arise is when looking at $S^r_d(\Delta)$.   Example~\ref{BilleraRoseMatrix} described a process to change a partition $\Delta$ into another partition $\hat \Delta$ with the property that $S^r(\hat \Delta)_d \cong S^r_d(\Delta)$ as modules. Taking duals, the spline space $S^r(G_{\hat \Delta})_d$ is related both to the quotient ring $S^r(G_{\hat \Delta})/\langle x_1, \ldots, x_n \rangle^d S^r(G_{\hat \Delta})$ and to the ring of splines over $G_{\hat \Delta}$ with coefficients in the quotient of the polynomial ring by $\langle x_1, \ldots, x_n \rangle^d$.  This quotient is an algebraic formalization of ``just ignore terms of degree at least $d$" but it endows $S^r_d(\Delta)$ with a ring structure.  We can ask when all of these maps are actually isomorphisms, which is related to questions about when the spline module is free.  We can also ask for the ring structure.

Passing from equivariant cohomology to ordinary cohomology is the same as taking the quotient by $\langle x_1, \ldots, x_n \rangle S^r(G)$.  This is different from restricting degree, though not unrelated.  For instance, when $X$ is the projective line $\mathbb{P}^1$ then its corresponding graph $G_X$ consists of a single edge, whose the label we can take to be $x \in \mathbb{R}[x]$.  Thus the elements in Figure \ref{Basis Spline Eqvt vs non Eqvt} generate $S^0(G_X)$ as a module over $\mathbb{R}[x]$.

\begin{figure}[h]\label{Basis Spline Eqvt vs non Eqvt}
\hspace{1cm}
\begin{picture}(60,20)(0,0)
\put(-30,10){$\sigma_1=$}
\put(0,0){\circle*{3}}
\put(0,20){\circle*{3}}
\put(0,0){\line(0,1){20}}
\put(3,0){$1$}
\put(3,20){$1$}
\end{picture}\hspace{1cm}
\begin{picture}(30,20)(0,0)
\put(-30,10){$\sigma_2=$}
\put(0,0){\circle*{3}}
\put(0,20){\circle*{3}}
\put(0,0){\line(0,1){20}}
\put(3,0){$0$}
\put(3,20){$x$}
\end{picture}
\caption{$\mathbb{R}[x]$-basis for the spline module $S^0(G)$ on a graph with one edge}
\end{figure}
This means every element in $S^0(G_X)$ has the form $a(x)\sigma_1+b(x)\sigma_2$ for $a(x),b(x)\in\mathbb{R}[x]$. The non-equivariant cohomology is obtained by evaluating the polynomial coefficients $a(x)$ and $b(x)$ at $x=0$ but taking the non-equivariant images $\widetilde{\sigma_1}$ and $\widetilde{\sigma_2}$, which correspond to the Poincar\'e dual of the class of the whole variety and of a point, respectively. Therefore $\widetilde{\sigma_2}$ is not zero, even though its only entries are $x$ and zero.  We would like efficient algorithms to pass from the splines that arise in equivariant cohomology to the ordinary cohomology.  More questions with a similar equivariant versus non-equivariant flavor can be found in, e.g., \cite{RYY}.

\subsubsection{Questions about representations on splines}

Section 3 described a way to construct group actions on splines, both as modules and as vector spaces, since the group actions preserve the degree of each spline. Many questions remain about these representations, especially when $r > 0$.

The first question raised in that section is: which partitions $\Delta$ or edge-labeled graphs give rise to ``interesting" edge-labeled automorphisms?  As a first pass,  ``interesting" could mean that the graph automorphism is not the identity.  However, another interesting graph automorphism would be one that preserves a local symmetry, for instance an instance of an interior vertex like that of $\Delta_1$ in Figure~\ref{3planarTriangulations} of Example~\ref{ex:first}.  This could provide one plan of attack to prove that Schumaker's lower bound formula gives $S^1_3(\Delta)$ for planar triangulations.

The examples in Section 3 also showed  strategies to decompose a representation on splines into subrepresentations. Given a representation of a group $H$, a classic question asks how to decompose the representation into particular kinds of subrepresentations, e.g. irreducibles (the smallest building blocks), representations arising from subgroups of $H$, etc.  We can ask this question about a representation on splines for any particular edge-labeled graph.

For instance, the last few years have seen particularly exciting developments for representations on splines for graphs $G_h$ that represent a family of varieties called \textit{(regular semisimple) Hessenberg varieties}. The \textit{Stanley-Stembridge conjecture} is a powerful combinatorial conjecture about immanants, which are certain functions on square matrices that generalize determinants.

The conjecture appears in evolving forms in \cite{stanStem, stan95}, as well as references in those papers to Stembridge’s earlier work on immanants, and Haiman’s proof that immanants of Jacobi-Trudi matrices can be expanded in terms of Schur functions. Shareshian and Wachs subsequently suggested that the Stanley-Stembridge conjecture could be reformulated as a statement about precisely how the representations on splines for the Hessenberg graphs $G_h$ decompose into irreducibles \cite{ShaWac16}.  Brosnan and Chow \cite{BroCho18}, and Guay-Paquet \cite{Gua2} proved the Shareshian-Wachs conjecture using very different methods.  The representations have since been shown to have deep connections to other combinatorial questions about, e.g., LLT polynomials. For history and evolution of the problem, see the survey by Shareshian and Wachs \cite{ShaWac2}.

 \subsubsection{Geometric interpretation of $r > 0$}
We close with a question for algebraists and geometers that comes from analysis.  The spline spaces $S^r(\Delta)$ for $r>0$ are very natural from the perspective of differentiable piecewise-polynomial functions.  Algebraically, it is also very natural to raise an edge-label to higher power, and it's clear that  we have a chain of inclusions
\[
S^0(\Delta) \subseteq S^1(\Delta) \subseteq S^2(\Delta) \cdots
\]
from which we can form the associated graded object. However natural this construction, there is currently no known geometric interpretation of this object. For example, is it connected to the process of taking blow-ups?  One possible first step towards this question is in \cite{SDCG}, which uses barycentric coordinates to translate the smoothness conditions to the Stanley-Reisner ring.  We would like to know more, even just for special cases like when $\Delta$ is the simplicial complex corresponding to projective space $\mathbb{P}^n$.

\vskip .1in
\noindent {\bf Acknowledgments}:  The authors thank the Istituto Nazionale di Alta Matematica (INdAM) for supporting the conference 
``Approximation Theory and Numerical Analysis meet Algebra, Geometry, Topology", which led to this work, the staff at the 
Palazzone di Cortona where the conference took place, and our fellow participants for providing a stimulating and enjoyable experience. Computations were performed using the computer algebra package {\tt Macaulay2} \cite{danmike}. 

ML acknowledges the PRIN2017 CUP E8419000480006, the Fondi di Ricerca Scientifica di Ateneo 2021 CUP E853C22001680005, and the MUR Excellence Department Project awarded to the Department of Mathematics, University of Rome Tor Vergata, CUP E83C23000330006.  HS acknowledges support from NSF 2006410 and the Rosemary Kopel Brown endowment at Auburn. JT acknowledges support from the NSF under grant number DMS-2054513 and from the AWM under their Mathematical Endeavors Revitalization Program.
\renewcommand{\baselinestretch}{1.0}
\small\normalsize 

\bibliographystyle{amsalpha}

\end{document}